\documentclass[reqno]{amsart2000}

\theoremstyle{plain}
\newtheorem{theorem}{Theorem}[section]
\newtheorem*{nntheorem}{Theorem}
\newtheorem{proposition}[theorem]{Proposition}

\newtheorem*{nnlemma}{Lemma}

\newtheorem*{corollary}{Corollary}

\theoremstyle{definition}
\newtheorem{definition}{Definition}[section]

\newtheorem*{nnexample}{Example}
\newtheorem*{remark}{Remark}
\newtheorem*{notation}{Notation}

\numberwithin{equation}{section}

\newcommand{\Catname}[1]{\mathbf{#1}}

\newcommand{\Complex}{\mathbb{C}}

\newcommand{\Field}{\mathbb{K}}

\usepackage[all]{xy}

\begin{document}

\title[Homotopy Quantum Field Theories \ldots]{Homotopy Quantum Field Theories and the Homotopy Cobordism Category in Dimension $1+1$.}

\author{Gon\c{c}alo Rodrigues}

\address{Centro de Matem\'{a}tica Aplicada\\
Departamento de Matem\'{a}tica\\
Instituto Superior T\'{e}cnico\\
Av. Rovisco Pais\\
1096, Lisboa\\
Portugal}

\email{grodr@math.ist.utl.pt}

\thanks{This work was done with the support of the grant PRAXIS
XXI/BD/17226/98 from \emph{Funda\c{c}\~{a}o para a Ci\^{e}ncia e Tecnologia}.}

\thanks{This work was also supported by the programme \emph{Programa
Operacional ``Ci\^{e}ncia, Tecnologia, Inova\c{c}\~{a}o''} (POCTI)
of the \emph{Funda\c{c}\~{a}o para a Ci\^{e}ncia e Tecnologia}
(FCT), cofinanced by the European Community fund FEDER}

\begin{abstract}
We define Homotopy quantum field theories (HQFT) as Topological
quantum field theories (TQFT) for manifolds endowed with extra
structure in the form of a map into some background space $X$. We
also build the category of homotopy cobordisms
$\Catname{HCobord}(n,X)$ such that an HQFT is a functor from this
category into a category of linear spaces. We then derive some
very general properties of $\Catname{HCobord}(n,X)$, including the
fact that it only depends on the $(n+1)$-homotopy type of $X$. We
also prove that an HQFT with target space $X$ and in dimension
$n+1$ implies the existence of geometrical structures in $X$; in
particular, flat gerbes make their appearance. We give a complete
characterization of $\Catname{HCobord}(n,X)$ for $n=1$ (or the
$1+1$ case) and $X$ the Eilenberg-Maclane space $K(G,2)$. In the
final section we derive state sum models for these HQFT's.
\end{abstract}
\maketitle

\section*{Introduction.}

The subject of Topological Quantum Field Theories is well established by now. Very early, homotopy theory methods were used to build examples of TQFT's. We can cite \cite{DW} and the works of D.~Yetter, \cite{DY} and \cite{DY1} (see also their reformulation and generalization by T.~Porter, \cite{TP} and \cite{TP1}). Essentially, what all these theories do is fix a background space $X$ and compute a weighted sum over homotopy classes of maps $f:M\longrightarrow X$ for a closed manifold $M$.

In all these cases the homotopy information is an auxiliary to build TQFT's. In \cite{TH1}, V.~Turaev shifted the status of the homotopy information by introducing HQFT's as an enrichment of Topological quantum field theories (TQFT) for $n$-manifolds and $(n+1)$-cobordisms endowed with extra structure in the form of a continuous map into some target space $X$.

One of the important theorems in \cite{TH1} was that the HQFT's only depended on the
$n$-homotopy type of $X$. Around the same time appeared the paper
\cite{BT} which discussed HQFT's in dimension $1+1$, but for a
simply connected target space, and therefore not covered by V.~Turaev's definition.

Here we present a new and broader definition of a Homotopy quantum
field theory in such a way that it only depends on the
$(n+1)$-homotopy type of $X$, enlarging Turaev's definition and
also covering the special HQFT's in \cite{BT}. This certainly was
one of our motivations for starting this work but it was not the
only one. If we consider the ``physical'' situation of $3+1$
dimensions then, once again quoting Turaev's theorem, the HQFT's
would only depend on the $3$-homotopy type of the target space. In
other words, we would be unable to detect $\pi_{4}(X)$. But
various models appearing in physics are related to $\pi_{4}(X)$ in
some way. Here we just mention E.~Witten's work reinterpreting the
Donaldson-Floer theory as a field theory, more precisely the
Donaldson polynomials, obtained through intersection theory in the
Instanton moduli space, were obtained, at least formally, as
correlation functions for a quantum field theory with a twisted
supersymmetric topological lagrangian (see \cite{WIT}).

Also on the physics side, in \cite{BC} J.~Barrett and L.~Crane
gave a model for euclidean quantum gravity, by starting with a
state sum model for a TQFT using the monoidal category of
representations of $SO(4)$ and then imposing a constraint in the
form of restricting the sum to be only over a specific
subcategory. Still, this is only a model for pure gravity, one
still has to insert matter. Here we take the view that just as
TQFT's can be considered as a first approximation to full-blown
quantum gravity, HQFT's are a first approximation to gravity
coupled with matter. Thus, studying abstract HQFT's is a first
step in trying to insert matter in state sum models.

\smallskip
The contents of the paper are as follows. In the first section we
give our definition of what is a Homotopy quantum field theory and
derive some very general properties. We also ask the question if
there is some kind of category of homotopy cobordisms
$\Catname{HCobord}(n,X)$ such that an HQFT is a functor from this
category into the category of linear spaces $\Catname{Vect}$. We
answer this question in the affirmative by providing a rigorous
construction of it -- due to the technicalities involved we have
relegated the construction to an appendix. This also has the very
practical consequence that our computations will be made inside
this category with no reference to any HQFT. We also prove that
this construction is functorial in $X$ and that
$\Catname{HCobord}(n,X)$ (and therefore HQFT's) only depends on
the $(n+1)$-homotopy type of $X$. Finally, we show how an HQFT
determines certain geometric structures, gerbes in particular, in
the target space.

In the second section we solve the $n=1$ (or $1+1$) case where the
target space $X$ is an Eilenberg-Maclane space $K(G,2)$, by giving
a complete characterization of $\Catname{HCobord}(1,X)$. In this
case the solution is particularly elegant, in that it is in terms
of a universal property.

In the third section we give a complete derivation of the state
sum models for these HQFT's. If $A$ is an associative semi-simple
algebra then we can enrich the corresponding (1+1)-TQFT to an HQFT
with target space $K(Z(A)^{\ast},2)$, where $Z(A)^{\ast}$ is the
group of invertible elements of the center of $A$. A homomorphism
$G\longrightarrow Z(A)^{\ast}$ allows the target space to be
``reduced'' to $K(G,2)$. All homotopical state sums are obtained
in this way.

\section{Homotopy quantum field theories: definition and basic
properties.}\label{section:defsandgenprops}

Before giving the full definition of a Homotopy Quantum Field
Theory (HQFT for short) let us settle some points on notation and
terminology.

\smallskip
When we use the word manifold we always mean a compact, oriented,
differentiable manifold. In particular, it admits at least one
structure of a finite CW-complex with cells up to the dimension of
the manifold.

\smallskip
Fix a path-connected space $X$ with a fixed base point $*$.

\begin{definition}
An $X$-manifold is a pair $(M,g)$ where $M$ is a closed manifold
with a choice of a base point $m_{i}$ for each connected component
$M_{i}$ of $M$. The discrete subspace of base points of $M$ will
be denoted by $\Catname{B}_{M}$. The letter $g$ stands for a
continuous map $M\longrightarrow X$, the characteristic map, such
that
\begin{equation*}
  g(m_{i})=\ast
\end{equation*}
that is, all the base points of $M$ are sent into the base point
$\ast$ of $X$.

An $X$-diffeomorphism $\psi :(M,g)\longrightarrow (N,h)$ between
$X$-manifolds is a diffeomorphism $\psi :M\longrightarrow N$
preserving the orientation and taking base points into base
points,
\begin{equation*}
  \psi(m_{i})=n_{j}
\end{equation*}
and such that the diagram in figure \ref{diag:commutefordiff} is
commutative\footnote{Throughout the paper we will always use the
diagrammatic order when writing a composition of morphisms.}.
\begin{figure}[htbp]
$\vcenter{\xymatrix{
  M \ar[rr]^{\psi} \ar[dr]_{f}
                &  &    N \ar[dl]^{g}    \\
                & X                 }
                }$
  \caption{}
  \label{diag:commutefordiff}
\end{figure}
\end{definition}

Denote by $\Catname{Diff}(n,X)$ the category of $n$-dimensional
$X$-manifolds as objects and $X$-diffeomorphisms as morphisms. We
define the disjoint union of $X$-manifolds by
\begin{equation*}
  (M,g)\amalg(N,h)=(M\amalg N,g\amalg h)
\end{equation*}
with the obvious definition for the map $g\amalg h:M\amalg
N\longrightarrow X$. With the operation of disjoint union
$\Catname{Diff}(n,X)$ is a symmetric monoidal category, where the
unit is the empty $X$-manifold $\emptyset$.

\begin{remark}
We warn the reader that we \emph{are not} identifying
$(M\amalg\emptyset,g)$ (or $(\emptyset\amalg M,g)$) with $(M,g)$.
Instead, they are naturally isomorphic via the obvious
$X$-diffeomorphism $l_{(M,g)}:(M\amalg\emptyset,g)\longrightarrow
(M,g)$ (or $r_{(M,g)}$).
\end{remark}

\begin{definition}
A cobordism $W:M\longrightarrow N$ between closed manifolds $M$
and $N$ is a manifold $W$ such that
\begin{equation*}
\partial W=M\amalg N
\end{equation*}
where $M$ and $N$ have chosen orientations via their normal
bundles, and such that the orientation induced from $W$ agrees
with that of $N$ and disagrees with the one on $M$. This state of
affairs is described symbolically by $\partial W= -M\amalg N$. We
will also use the notations $\partial_{-}W=M$ and
$\partial_{+}W=N$ to denote the incoming and the outgoing
components of the boundary.

An $X$-cobordism $(W,F):(M,g)\longrightarrow (N,h)$ is a cobordism
$W:M\longrightarrow N$ with a homotopy class of maps
$F:W\longrightarrow X$ relative to the boundary, such that
$F|_{M}=g$ and $F|_{N}=h$. We note that throughout the paper we
will make no notational distinction between the homotopy class $F$
and any of its representatives.

Finally, an $X$-diffeomorphism between $X$-cobordisms $(W,F)$ and
$(V,G)$ is a diffeomorphism $\Psi:W\longrightarrow V$ such that
\begin{equation*}
\begin{split}
\Psi(\partial_{+}W) &=\partial_{+}V \\
\Psi(\partial_{-}W) &=\partial_{-}V
\end{split}
\end{equation*}
and $F=\Psi G$, where equality is understood as equality of
homotopy classes of maps.
\end{definition}

\smallskip
Let us examine these definitions a little closer. Let $W$ be a
cobordism with no boundary, then an $X$-structure on it is just a
homotopy class of maps $W\longrightarrow X$. Now if we vary $X$ we
get different types of ``fields'' in $W$. For example, if $X$ is
the classifying space $BG$ (or $K(G,1)$ when $G$ is discrete) for
some group $G$ then we get an isomorphism class of $G$-bundles
over $W$. On the other hand if we put $X$ to be some
Eilenberg-Maclane space $K(G,n)$, where $G$ is an abelian group
and $n>1$ then the isomorphism
\begin{equation*}
  [W,K(G,n)]\cong H^{n}(W;G)
\end{equation*}
tells us that the ``fields'' are classes in the cohomology of $W$
with values in $G$. Finally if $X$ is just the $1$-point space
$\ast$ then we revert to the original set-up of TQFT's.

In the usual differential geometric description, quantum fields
are connections in bundles. More precisely, the interaction
fields are connections on principal bundles and the matter fields
are sections of the associated vector bundles. A direct analogy can be drawn in the HQFT case. The Gauge group
role is played by the background space X (take $K(G,1)$, for
example), the interaction fields are the maps into the background
space (in the $K(G,1)$ case they are exactly $G$-Principal
bundles) and the matter fields appear by taking representations of the category
$\Catname{HCobord}(n,X)$ in $\Catname{Vect}$ (in the $K(G,1)$
case, taking the associated bundle by some representation.)

There is another point in saying that we have something like
matter surfacing. In the case where the background space is
$K(G,n)$, the maps are just cohomology classes Poincar\'{e} dual to 
homology classes, which are localizable, similarly to local matter fields. Since every (simple) background space can be ``synthesized'' in $K(G,n)$ spaces via Postnikov towers, the analogy extends to these spaces as well.

\smallskip
There are two basic operations we can perform with cobordisms that
extend to $X$-cobordisms: gluing and disjoint union. Given two
$X$-cobordisms $(W,F):(M,g)\longrightarrow (N,h)$ and
$(V,G):(N',h')\longrightarrow (P,j)$, and an $X$-diffeomorphism
$\psi:(N,h)\longrightarrow(N',h')$, we have a composition
cobordism $(W\amalg_{\psi}V,F\cdot G ):(M,g)\longrightarrow (P,j)$
where $V\amalg_{\psi}W$ is the gluing of $V$ and $W$ along their
boundary components identified by $\psi$ and $F\cdot G$ is defined
by
\begin{equation*}
F\cdot G(x)=\begin{cases} F(x)& \text{if $x\in W$}\\
G(x)& \text{if $x\in V$}
\end{cases}
\end{equation*}

We leave to the reader the easy task of checking that this is well
defined. If $\psi$ is the identity we denote the gluing by
$(W,F)\circ(V,G)$.

\smallskip
For each $X$-manifold $(M,g)$ there is also an $X$-cobordism
$(I\times M,1_{g}):(M,g)\longrightarrow(M,g)$ with
$1_{g}(t,x)=g(x)$ and where $I$ is the unit interval. This
cobordism will be called the identity cobordism and denoted by
$1_{(M,g)}$.

\smallskip
As for the disjoint union $(W,F)\amalg(V,G)$ of $X$-cobordisms it
is defined in the same way as the disjoint union of $X$-manifolds.

\smallskip
After this introduction we can state the definition of a homotopy
quantum field theory.

\begin{definition} \label{def:hft}
Fix a commutative ring $R$ and denote by $\Catname{Mod}(R)$ the
category of free finitely generated $R$-modules. Usually, the ring
$R$ is the complex field $\Complex$, and then
$\Catname{Mod}(\Complex)$ is just the category
$\Catname{Vect}(\Complex)$ of complex linear spaces.

An $n$-dimensional Homotopy quantum field theory $\tau$ with
target space $X$ is an assignment of a finitely generated free
$R$-module $\tau(M,g)$ to every $X$-manifold $(M,g)$ of dimension
$n$, a linear isomorphism
$\tau\psi:\tau(M,g)\longrightarrow\tau(N,h)$ to every
$X$-diffeomorphism $\psi:(M,g)\longrightarrow(N,h)$, and a linear
map $\tau(W,F):\tau(M,g)\longrightarrow\tau(N,h)$ to every
$X$-cobordism $(W,F):(M,g)\longrightarrow(N,h)$. This data
satisfies the following axioms
\begin{enumerate}
\item \label{axm:ax1}
$\tau$ is functorial in $\Catname{Diff}(n,X)$. This means that for
two $X$-diffeomorphisms $\psi:(M,g)\longrightarrow(N,h)$ and
$\phi:(N,h)\longrightarrow(P,j)$ we have
\begin{equation*}
  \tau(\psi\phi)=\tau(\psi)\tau(\phi)
\end{equation*}
and if $1_{(M,g)}$ is the identity $X$-diffeomorphism on $(M,g)$
then $\tau(1_{(M,g)})=1_{\tau(M,g)}$.

\item \label{axm:ax2}
$\tau$ is also symmetric monoidal in $\Catname{Diff}(n,X)$. This
means that there are natural isomorphisms
\begin{equation*}
  c_{(M,g),(N,h)}:\tau((M,g)\amalg(N,h))\cong\tau(M,g)\otimes\tau(N,h)
\end{equation*}
and an isomorphism $u:\tau(\emptyset)\cong R$ that satisfy the
usual axioms for a symmetric monoidal functor (see \cite{AJRS} for
the actual diagrams).

\item \label{axm:ax3}
For $X$-cobordisms $(W,F):(M,g)\longrightarrow(N,h)$ and
$(V,G):(N',h')\longrightarrow(P,j)$ glued along
$\psi:(N,h)\longrightarrow(N',h')$ we have
\begin{equation*}
  \tau((W,F)\amalg_{\psi}(V,G))=\tau(W,F)\tau(\psi)\tau(V,G)
\end{equation*}

\item \label{axm:ax4}
For the identity $X$-cobordism $1_{(M,g)}=(I\times M,1_{g})$ we
have
\begin{equation*}
  \tau1_{(M,g)}=1_{\tau(M,g)}
\end{equation*}

\item \label{axm:ax5}
For $X$-cobordisms $(W,F):(M,g)\longrightarrow(N,h)$, $(V,G):(M',g')\longrightarrow(N',h')$ and $(P,j):\emptyset\longrightarrow\emptyset$, we have the commutative
diagrams in figure \ref{diag:natforcobord}.
\begin{figure}[hbtp]
  $\vcenter{\xymatrix{
    \tau((M,g)\amalg(M',g')) \ar[d]_{\tau((W,F)\amalg(V,G))} \ar[r]^{c} & \tau(M,g)\otimes\tau(M',g')  \ar[d]^{\tau(W,F)\otimes\tau(V,G)} \\
    \tau((N,h)\amalg(N',h')) \ar[r]^{c} & \tau(N,h)\otimes\tau(N',h')   }}$
  $\vcenter{\xymatrix{
  \tau\emptyset \ar[r]^{u} \ar[d]_{\tau(P,j)} & R  \\
  \tau\emptyset \ar[ur]_{u}  & }
  }$
  \caption{}
  \label{diag:natforcobord}
\end{figure}

\item \label{axm:ax6}
For every $X$-diffeomorphism $\Psi:(W,F)\longrightarrow(V,G)$
between $X$-cobordisms we have the commutative diagram in figure
\ref{diag:xdiffrel}.
\begin{figure}[hbtp]
  $\vcenter{\xymatrix@C+15pt{
    \tau(\partial_{-}(W,F)) \ar[d]_{\tau(W,F)} \ar[r]^{\tau(\Psi|_{\partial_{-}W})} & \tau(\partial_{-}(V,G)) \ar[d]^{\tau(V,G)} \\
    \tau(\partial_{+}(W,F)) \ar[r]^{\tau(\Psi|_{\partial_{+}W})} & \tau(\partial_{+}(V,G))   }}$
  \caption{}
  \label{diag:xdiffrel}
\end{figure}
\end{enumerate}
\end{definition}

As we have already mentioned in the introduction there are some
differences between our definition of HQFT and that of V.~Turaev in
\cite{TH1}. In axiom \ref{axm:ax2} we have added the structural
isomorphisms and in axiom \ref{axm:ax5} we state that these
isomorphisms are natural \emph{also} for cobordisms. This
naturality condition will be essential for the successful
construction of the category of ``homotopy cobordisms.''

A more substantial difference is what we take as identity
cobordisms in axiom \ref{axm:ax4}. Our definition is much stricter
and implies that there are many more cylinder isomorphisms in our
theory. While in \cite{TH1} for every cylinder $I\times M$ the
homotopical information on it is completely determined by the
information on the boundary, in our definition this specification
is not enough (see propositions \ref{thm:modinv} and
\ref{thm:isodescription} below). And it is this fact that
ultimately justifies the theorem that our HQFT's will depend on
the $(n+1)$-homotopy type of $X$, instead of the $n$-homotopy type
as in \cite{TH1}. It is also this fact that implies that an HQFT
determines a (flat) gerbe in the target space $X$.

\smallskip
Back to the definitions. There is also a notion of map between two
HQFT's.

\begin{definition}
Let $\tau$ and $\rho$ be two HQFT's in dimension $n+1$ and with
target space $X$, then a map $\theta:\tau\longrightarrow\rho$ is a
family of maps $\theta_{(M,g)}:\tau(M,g)\longrightarrow\rho(M,g)$
indexed by $X$-manifolds $(M,g)$, such that for every
$X$-diffeomorphism $\psi:(M,g)\longrightarrow(N,h)$ and every
$X$-cobordism $(W,F):(M,g)\longrightarrow(N,h)$ we have the
commutative diagrams in figure \ref{diag:natfortransf} (the
naturality conditions) and also the commutative diagrams in figure
\ref{diag:nattensor}.
\begin{figure}[hbtp]
  $\vcenter{\xymatrix{
    \tau(M,g) \ar[d]_{\tau\psi} \ar[r]^{\theta_{(M,g)}} & \rho(M,g) \ar[d]^{\rho\psi} \\
    \tau(N,h) \ar[r]^{\theta_{(N,h)}} & \rho(N,h)   }}$
  $\vcenter{\xymatrix{
    \tau(M,g) \ar[d]_{\tau(W,F)} \ar[r]^{\theta_{(M,g)}} & \rho(N,h) \ar[d]^{\rho(W,F)} \\
    \tau(N,h) \ar[r]^{\theta_{(N,h)}} & \rho(N,h)   }}$
  \caption{}
  \label{diag:natfortransf}
\end{figure}

\begin{figure}[hbtp]
  $\vcenter{\xymatrix{
    \tau\left((M,g)\amalg(N,h)\right) \ar[d]_{\theta_{(M,g)\amalg(N,h)}} \ar[r]^{c} & \tau(M,g)\otimes\tau(N,h) \ar[d]^{\theta_{(M,g)}\otimes\theta_{(N,h)}} \\
    \rho\left((M,g)\amalg(N,h)\right) \ar[r]_{c} & \rho_{(M,g)}\otimes\rho_{(N,h)}   }}$
  $\vcenter{\xymatrix{
  \tau\emptyset \ar[d]_{\theta_{\emptyset}} \ar[dr]^{u_{\tau}}        \\
  \rho\emptyset \ar[r]_{u_{\rho}}  & R              }
  }$
  \caption{}
  \label{diag:nattensor}
\end{figure}
\end{definition}

It is easy to check that there is a category $\Catname{HQFT}(n,X)$
whose objects are $n$ dimensional HQFT's with target space $X$ and
the morphisms are the maps between them as defined above.

\smallskip
If we look at the definition of an HQFT we see that it definitely
has a functorial look. For example, axioms \ref{axm:ax1} and
\ref{axm:ax2} say explicitly that $\tau$ is a symmetric monoidal
functor on $\Catname{Diff}(n,X)$. So the question is, does there
exist a symmetric monoidal category $\Catname{HCobord}(n,X)$ such
that we have an isomorphism\footnote{An equivalence of categories
would suffice, but from the actual construction it is an
isomorphism that drops out.}
\begin{equation*}
  \Catname{HQFT}(n,X)\cong[\Catname{HCobord}(n,X),\Catname{Mod}(R)]
\end{equation*}
where $[\mathcal{A},\mathcal{B}]$ is the category of symmetric
monoidal functors $\mathcal{A}\longrightarrow\mathcal{B}$? The
answer is yes. In order not to interrupt the flow of the paper we
have relegated the construction to the appendix and here we
content ourselves with making some comments on the rationale
behind the category.

\smallskip
The objects of $\Catname{HCobord}(n,X)$ are clearly $X$-manifolds
but since an HQFT $\tau$ assigns linear maps to both
$X$-diffeomorphisms and $X$-cobordisms the morphisms in
$\Catname{HCobord}(n,X)$ have to be a mixture of both.

Formally, there are two types of axioms. Those involving the
structural isomorphisms, like axiom \ref{axm:ax2}, that express
the fact that the assignments are functorial, and those that only
involve applications of $\tau$, like axioms \ref{axm:ax3} and
\ref{axm:ax6}. Let us look more attentively at one of them, at the
equation
\begin{equation*}
  \tau((W,F)\amalg_{\psi}(V,G))=\tau(W,F)\tau(\psi)\tau(V,G)
\end{equation*}
in axiom \ref{axm:ax3}. We have $\tau$ applied on both sides and
on all terms, so it must express a general relation in our
would-be category $\Catname{HCobord}(n,X)$ of ``homotopy
cobordisms''. Let us peel off the $\tau$ from the equation to get
\begin{equation*}
  (W,F)\amalg_{\psi}(V,G)=(W,F)\psi(V,G)
\end{equation*}

On the left hand side we have an $X$-cobordism and on the right we
have a formal composition of two $X$-cobordisms and an
$X$-diffeomorphism. This hints at the fact that morphisms in
$\Catname{HCobord}(n,X)$ are alternating strings of
$X$-diffeomorphisms and $X$-cobordisms. Axioms for an HQFT (like
\ref{axm:ax3} and \ref{axm:ax6}) express universal relations
between strings of these morphisms. Up to some technical details
dealing with the lack of associativity for gluing cobordisms,
this \emph{is} the right image to have in mind and it will be
implicitly used in the proofs below where we consider the cases
of  $X$-cobordisms and $X$-diffeomorphisms one at a time.

Although the construction of the category $\Catname{HCobord}(n,X)$
is left for the appendix, we need to spell out the most important
relations holding in this category.

\begin{notation}
In this list, and in the rest of the paper with the exception of
the appendix, we make an abuse of notation and identify an
$X$-cobordism $(W,F)$ and an $X$-diffeomorphism $\psi$ with their
images in the category $\Catname{HCobord}(n,X)$.
\end{notation}

They are specified in the next list.

\begin{enumerate}
\item \label{reln:xdiffrel}
If $\Psi:(W,F)\longrightarrow(V,G)$ is an $X$-diffeomorphism
between $X$-cobordisms then we have
\begin{equation*}
  (W,F)\Psi|_{\partial_{+}W}=\Psi|_{\partial_{-}W}(V,G)
\end{equation*}

\item \label{reln:gluewithf}
If we have two $X$-cobordisms $(W,F):(M,g)\longrightarrow(N,h)$
and $(V,G):(N',h')\longrightarrow(P,j)$ and an $X$-diffeomorphism
$\psi:(N,h)\longrightarrow(N',h')$ we have
\begin{equation*}
  (W,F)\amalg_{\psi}(V,G)=(W,F)\psi(V,G)
\end{equation*}

\item \label{reln:ident}
For $X$-cobordisms $(W,F):(M,g)\longrightarrow(N,h)$ and
$(V,G):(N,h)\longrightarrow(P,j)$ we have
\begin{equation*}
\begin{split}
  (W,F)\circ(I\times N,1_{h}) &=(W,F) \\
  (I\times N,1_{h})\circ(V,G) &=(V,G)
\end{split}
\end{equation*}
and $(I\times N,1_{h})$ equals the identity $X$-diffeomorphism
$1_{(N,h)}$.
\end{enumerate}

\smallskip
After this short introduction to $\Catname{HCobord}(n,X)$ we start
establishing some of its properties. Since $X$-diffeomorphisms are
among its morphisms one could believe that such a category would
be very big. The next proposition shows that such a belief is
incorrect.

\begin{proposition} \label{thm:isotinv}
If $\psi$ and $\phi$ are isotopic $X$-diffeomorphisms
$(M,g)\longrightarrow(N,h)$ then they are equal in
$\Catname{HCobord}(n,X)$.
\end{proposition}
\begin{proof}
Remember that $\psi$ and $\phi$ are isotopic iff there is a
(differentiable) homotopy $H:I\times M\longrightarrow N$ between
them such that for fixed $t\in I$, $H_{t}(x)$ is an
$X$-diffeomorphism $(M,g)\longrightarrow(N,h)$. Consider then the
identity cobordisms $(I\times M,1_{g})$ and $(I\times N,1_{h})$.
The map $\pi_{I}\times H:I\times M\longrightarrow I\times N$,
where $\pi_{I}$ is the projection $I\times M\longrightarrow I$, is
clearly an $X$-diffeomorphism. Since when restricted to the
boundary it is the disjoint union of $\psi$ and $\phi$, relation
\ref{reln:xdiffrel} gives the desired result.
\end{proof}

This proposition is one of the hallmarks of a topological theory.
The action of the full diffeomorphism group is filtered and only
the mapping class group (or the $X$ version of it) survives.

The next proposition shows that, up to isomorphism, $(M,g)$
depends only on the homotopy class of $g$. More precisely

\begin{proposition} \label{thm:modinv}
Let $(M,g)$ and $(M,h)$ be $X$-manifolds such that there is a free
homotopy between $g$ and $h$, then there is an isomorphism
$(M,g)\longrightarrow(M,h)$ in $\Catname{HCobord}(n,X)$.
\end{proposition}
\begin{proof}
Let $F:I \times M \longrightarrow X$ be an homotopy between $g$
and $h$. Then we have an $X$-cobordism $(I\times
M,F):(M,g)\longrightarrow(M,h)$. On the other hand if we define
$\overline{F}:I\times M\longrightarrow X$ by
\begin{equation*}
  \overline{F}(t,x)=F(1-t,x)
\end{equation*}
we have another cobordism $(I\times
M,\overline{F}):(M,h)\longrightarrow(M,g)$. What we are going to
show is that they are inverse to each other. Of course we can
"shift" this cobordism to $([1,2]\times M,\overline{F})$ with
$\overline{F}$ now defined by $\overline{F}(t,x)=F(2-t,x)$, by the
diffeomorphism $(t,x)\longmapsto (1+t,x)$. Now we perform the
gluing of these two $X$-cobordisms to obtain $(I\times
M,F)\circ([1,2]\times M,\overline{F})=([0,2]\times
M,F\cdot\overline{F})$ and we will prove that it is the identity
by first deforming $F\cdot\overline{F}$ to a ''constant'' map and
then using a contracting diffeomorphism $[0,2]\longrightarrow I$
and relation \ref{reln:xdiffrel}.

A homotopy $1_{g}\longrightarrow F\cdot\overline{F}$ is a map
$H:I\times([0,2]\times M)\longrightarrow X$. Define $H$ by
\begin{equation*}
  H(s,(t,x))=
  \begin{cases}
    F(t,x) & \text{if $0\leq t\leq s$}, \\
    F(s,x) & \text{if $s\leq t\leq 2-s$}, \\
    \overline{F}(t,x) & \text{if $2-s\leq t\leq 2$}.
  \end{cases}
\end{equation*}

Note that $H(0,(t,x))=1_{g}(t,x)=g(x)$ and
$H(1,(t,x))=F\cdot\overline{F}(t,x)$. The fact that it is a
well-defined and continuous map and also a homotopy relative to
the boundary are also easy-to-establish facts that we leave to the
reader. This implies that $([0,2]\times
M,F\cdot\overline{F})=([0,2]\times M,1_{g})$. Now define the
diffeomorphism $\psi:[0,2]\longrightarrow I$ by
\begin{equation*}
  \psi(t)=\frac{1}{2}t
\end{equation*}

This extends to an $X$-diffeomorphism $\psi\times
1_{M}:([0,2]\times M,1_{g})\longrightarrow(I\times M,1_{g})$ which
is the identity on the boundary, so, by applying relation
\ref{reln:xdiffrel} we have finally $([0,2]\times
M,F\cdot\overline{F})=(I\times M,1_{g})$.

This gives $(I\times M,F)\circ(I\times
M,\overline{F})=1_{(M,g)}$. To get the equality in the other
direction first shift by $I\longrightarrow[-1,0]$ and rerun the
above proof with suitable changes.
\end{proof}


We can improve on proposition \ref{thm:modinv} a little bit. Recall that an (admissible) Morse function on a cobordism $W$ is a Morse function $W\longrightarrow I$ with no critical values on $\partial W$ and with
\begin{align*}
  f(\partial_{-}W) &= {0} \\
  f(\partial_{+}W) &= {1}
\end{align*}
if $\partial_{-}W$ ($\partial_{+}W$) is non-empty. Now, Call an $X$-cobordism $(W,F):(M,g)\longrightarrow (N,h)$ a cylinder iff $W$ is diffeomorphic to a cylinder $M\times I$. It is a cylinder
cobordism iff there is a Morse function on $W$ with no critical
values.

A cylinder cobordism is a cylinder, but a cylinder is not necessarily a cylinder cobordism -- just consider the handle $\emptyset\longrightarrow S^{1}\amalg S^{1}$ in two dimensions. The condition of $W$ having a Morse function with no critical values is precisely the one that guarantees that the source $M$ is diffeomorphic to the target $N$.

With these definitions we have


\begin{proposition} \label{thm:isodescription}
In the category $\Catname{HCobord}(n,X)$ we have
\begin{enumerate}

\item An $X$-diffeomorphism $\psi$ is equal to a cylinder cobordism.

\item A cylinder cobordism is an isomorphism.

\end{enumerate}
\end{proposition}
\begin{proof}
If $\psi$ is an $X$-diffeomorphism $(M,g)\longrightarrow(N,h)$, then by relation
\ref{reln:ident} and relation \ref{reln:gluewithf} we have
\begin{equation*}
  \psi=(I\times M,1_{g})\psi(I\times N,1_{h})=(I\times M\amalg_{\psi}I\times N,1_{g}\cdot 1_{h})
\end{equation*}
and the term on the right is clearly a cylinder cobordism.

For the second statement, note that if $W$ is a
cobordism with a Morse function $f$ that has no critical values
then by the cylinder recognition theorem (see the book \cite{MH},
page 153) there is a diffeomorphism $\Psi:M\times I\longrightarrow
W$ such that the triangle in figure \ref{diag:cylinderrecognition}
is commutative.
\begin{figure}[hbtp]
  $\vcenter{\xymatrix{
  M\times I \ar[dr]_{\pi_{I}} \ar[r]^{\Psi} & W \ar[d]^{f}  \\
                & I             }
  }$
  \caption{}
  \label{diag:cylinderrecognition}
\end{figure}

Since $\Psi$ is the identity in $M\times \{0\}$, relation
\ref{reln:xdiffrel} implies the equality
\begin{equation*}
  \Psi|_{M\times\{1\}}(W,F)=(M\times I,\Psi^{-1}F)
\end{equation*}

The right hand side is an isomorphism by proposition
\ref{thm:modinv}, and pushing $\Psi|_{M\times{1}}$, which is a
diffeomorphism, to the right side we have the result.
\end{proof}


\begin{remark}
It is not necessarily true that an isomorphism is a cylinder cobordism. It is possible that handle cancellations conspire to produce a cobordism with non-trivial topology and that at the same time has an inverse.
\end{remark}

\smallskip
In the above three propositions we have essentially proved that
$\Catname{HCobord}(n,X)$ is not so big after all, since everything
in sight only depends on the homotopy type. In the next
proposition we prove that $\Catname{HCobord}(n,X)$ has more
structure than just being symmetric monoidal, namely, it is a
category with duals. There are various definitions of what could
be a category with duals, with varying degrees of laxness built
into them, but the definition more suitable to our purposes is the
one in \cite{PFDY}.

In our discussion of duality we will make several statements whose
proof is omitted. To see those proofs, as well as other details on
the subject of duality in categories, the reader can look at, for
example, \cite{AJRS}, \cite{JBBW}, \cite{PFDY1} and \cite{TQ},
besides the reference cited above.

\begin{definition}
A monoidal category $(\mathcal{A},\otimes,\Field)$ has a duality
structure (or is a category with duals) iff for each object $a$
there is an object $a^{\ast}$, the dual object, and maps
$\eta_{a}:\Field\longrightarrow a\otimes a^{\ast}$,
$\varepsilon_{a}:a^{\ast}\otimes a\longrightarrow\Field$,
satisfying the following triangular identities\footnote{From now
on we suppress references to the structural isomorphisms of the
monoidal category. Maclane's coherence theorem justifies the
practice.}
\begin{equation} \label{eqn:triangident}
\begin{split}
  (\eta_{a}\otimes 1_{a})(1_{a}\otimes \varepsilon_{a}) &=1_{a} \\
  (1_{a^{\ast}}\otimes\eta_{a})(\varepsilon_{a}\otimes 1_{a^{\ast}}) &=1_{a^{\ast}}
\end{split}
\end{equation}
\end{definition}

In a category with duals, define the dual of a morphism
$f:a\longrightarrow b$, $f^{\ast}:b^{\ast}\longrightarrow
a^{\ast}$ by
\begin{equation*}
  f^{\ast}=(1_{b^{\ast}}\otimes\eta_{a})(1_{b^{\ast}}\otimes f\otimes
  1_{a^{\ast}})(\varepsilon_{b}\otimes 1_{b^{\ast}})
\end{equation*}

This definition of dual morphism lifts to a functor
$\mathcal{A}^{\ast}\longrightarrow\mathcal{A}$, where
$\mathcal{A}^{\ast}$ is the dual (or opposite) category of
$\mathcal{A}$. This functor is an equivalence of categories, and
if $\mathcal{A}$ is symmetric monoidal then there are natural
isomorphisms $(a\otimes b)^{\ast}\cong a^{\ast}\otimes b^{\ast}$,
$a^{\ast\ast}\cong a$ and $\Field^{\ast}\cong\Field$. In other
words, if we give to $\mathcal{A}^{\ast}$ the same monoidal
structure as $\mathcal{A}$ but where the structural isomorphisms
are the inverses of the ones of $\mathcal{A}$, then the dual
functor is a monoidal one.

\begin{definition}
A monoidal category $(\mathcal{A},\otimes,\Field)$ has strict
duals iff it has duals, the three structural isomorphisms
mentioned in the above paragraph are all identities and we have
the equality
\begin{equation} \label{eqn:strictdualeq}
  \eta_{a}=(\varepsilon_{a^{\ast}})^{\ast}
\end{equation}
plus
\begin{equation} \label{eqn:tensorprodrule}
  \eta_{a\otimes b}=\eta_{a}(1_{a}\otimes\eta_{b}\otimes 1_{a^{\ast}})(1_{a}\otimes 1_{b}\otimes \sigma_{a^{\ast},b^{\ast}})
\end{equation}
and a similar one for the pairing $\varepsilon$.
\end{definition}

\begin{remark}
The perspicacious reader will have noted that we have identified
$(a\otimes b)^{\ast}$ with $a^{\ast}\otimes b^{\ast}$ and
\emph{not} with $b^{\ast}\otimes a^{\ast}$ as is usually done. We
can get away with this because our monoidal categories are
symmetric and not just braided.
\end{remark}

\smallskip
The duality structure in $\Catname{HCobord}(n,X)$ is given by the
following: for $X$-manifolds $(M,g)$ we define $(M,g)^{\ast}$ by
$(-M,g)$, where $-M$ is the manifold $M$ but with the opposite
orientation. The maps $\eta_{(M,g)}$ and $\varepsilon_{(M,g)}$ are
defined by cylinders. More precisely, $\eta_{(M,g)}$ is the
cylinder $(I\times M,1_{g})$ viewed as an $X$-cobordism
$\emptyset\longrightarrow(M,g)\amalg(M,g)^{\ast}$, and
$\varepsilon_{(M,g)}$ is the cylinder $(I\times M,1_{g})$ viewed
as an $X$-cobordism
$(M,g)^{\ast}\amalg(M,g)\longrightarrow\emptyset$.

Before proving that this does in effect constitute a duality
structure for the category $\Catname{HCobord}(n,X)$, let me note
that there is another natural candidate for the functor part of
the duality. It is defined in the same way for objects $(M,g)$,
but for an $X$-cobordism $(W,F)$, its dual is the same
$X$-cobordism $(W,F)$ but now viewed as a map
$N^{\ast}\longrightarrow M^{\ast}$, and for an $X$-diffeomorphism
$\psi:M\longrightarrow N$ we have
\begin{equation*}
  \psi^{\ast}=\psi_{-}^{-1}
\end{equation*}
where $\psi_{-}$ is the diffeomorphism $-M\longrightarrow -N$
induced by $\psi$. Note that these ``duals'' automatically satisfy
the strictness condition. Below what we will prove is that both
dualities are in fact one and the same.

\begin{proposition} \label{thm:dualexist}
With the above definitions for a dual object and for the maps
$\eta$ and $\varepsilon$, $\Catname{HCobord}(n,X)$ is a monoidal
category with strict duals.
\end{proposition}
\begin{proof}
Note that $(\eta_{(M,g)}\amalg 1_{(M,g)})(1_{(M,g)}\amalg
\varepsilon_{(M,g)})$ is the gluing of four cylinders, each copy
with only trivial homotopy information on it. So we can first
''straighten'' out the gluing by a convenient diffeomorphism to
conclude the equality with $([0,4]\times M,1_{g})$ and then
squeeze $[0,4]$ to the unit interval $I$ to get the identity. The
other triangular identity is obtained in the same way.

Now let us prove that $(W,F)^{\ast}$ is just $(W,F)$ but as a
morphism $N^{\ast}\longrightarrow M^{\ast}$. By definition,
$(W,F)^{\ast}$ is obtained from $(W,F)$ by gluing four cylinders,
two at the beginning and the other two at the end. So, once again,
by ''straightening'' out the gluing, we get
\begin{equation*}
  (W,F)^{\ast}=([0,2]\times(-N),1_{h})(W,F)([0,2]\times(-M),1_{g})
\end{equation*}

Squashing $[0,2]$ to the unit interval $I$ we get the identity on
both sides, and the equality is proved.

Now let $\psi:(M,g)\longrightarrow(N,h)$ be an $X$-diffeomorphism,
then by proposition \ref{thm:isodescription} we have the equality
\begin{equation*}
  \psi=(I\times M\amalg_{\psi}I\times N,1_{g}\cdot 1_{h})
\end{equation*}

Applying the above calculation to the cobordism on the right we
have
\begin{equation*}
  \psi^{\ast}=(I\times M\amalg_{\psi}I\times N,1_{g}\cdot
  1_{h})^{\ast}=(I\times(-N)\amalg_{\psi_{-}}I\times(-M),1_{g}\cdot
  1_{h})
\end{equation*}

Now, a gluing made along $\psi$ is clearly diffeomorphic to a
gluing made along $\psi^{-1}$, so we have
\begin{equation*}
  \psi^{\ast}=(I\times(-N)\amalg_{\psi_{-}^{-1}}I\times(-M),1_{g}\cdot
  1_{h})=(I\times(-N),1_{g})\psi_{-}^{-1}(I\times(-M),1_{h})=\psi_{-}^{-1}
\end{equation*}

From these calculations it is clear that the duality is strict.
Equation \eqref{eqn:strictdualeq} is also obvious from the above
computation of duals. The equality \eqref{eqn:tensorprodrule}
follows from exactly the same type of arguments used up to now and
is left to the reader.
\end{proof}

There are some important consequences that we can extract from the
existence of a duality structure in $\Catname{HCobord}(n,X)$, or
on any monoidal category $\mathcal{A}$ for that matter. In the
first place, $\mathcal{A}$ is a closed category. There are natural
isomorphisms
\begin{equation}
\begin{split} \label{iso:adjointiso}
  \mathcal{A}(a\otimes b,c)\cong\mathcal{A}(b,a^{\ast}\otimes c) \quad&,\quad \mathcal{A}(a,b\otimes c)\cong\mathcal{A}(a\otimes c^{\ast},b) \\
  \mathcal{A}(a\otimes b,c)\cong\mathcal{A}(a,c\otimes b^{\ast})\quad&,\quad \mathcal{A}(a,b\otimes c)\cong\mathcal{A}(b^{\ast}\otimes a,c)
\end{split}
\end{equation}

We call the adjoint of $f$ the image of $f$ under any of the
isomorphisms \eqref{iso:adjointiso}. For example, the fourth
isomorphism is given explicitly by
\begin{equation*}
  f\longmapsto (1_{b^{\ast}}\otimes f)(\varepsilon_{b}\otimes 1_{c})
\end{equation*}
and there are similar formulas for the other isomorphisms. In
terms of our category $\Catname{HCobord}(n,X)$ it means that given
a cobordism with the boundary divided in a specific way, it can be
turned into another cobordism where the boundary has one more or
one less incoming or outgoing component by appropriately composing
with the duality pairings $\eta$ and $\varepsilon$. The
artificiality of viewing the same oriented cobordism in these
different ways is captured completely by the duality structure.

We also have the following two propositions that we quote from
\cite{AJRS}.

\begin{proposition}
Monoidal functors preserve duals, in particular, an HQFT $\tau$
takes the dual of an $X$-manifold $(M,g)$ to the dual of the
linear space $\tau(M,g)$ and the same for morphisms.
\end{proposition}

This proposition ``explains'' why duality did not enter in the
definition of an HQFT. If we view a symmetric monoidal category
with duals as the categorification of the notion of an abelian
group, then the above proposition is a categorification of the
simple fact that a morphism between two groups only has to
preserve the identity and the multiplication for it to
automatically preserve all the inverses. This paper is then a
study of the representation theory of some specific ``categorified
abelian groups.''

The other proposition is

\begin{proposition}
Every monoidal natural transformation between monoidal functors
defined on a category with duals is invertible, in particular,
every morphism between two HQFT's is in fact an isomorphism.
\end{proposition}

There is another piece of structure that we have not talked about:
the reversal of orientation for cobordisms. This structure goes
under the name of a $\ast$-structure of which we now give a formal
definition.

\begin{definition}
A $\ast$-structure in a monoidal category $\mathcal{A}$ is a
contravariant endofunctor $\dag:A\longrightarrow A$ which is the
identity on objects and satisfies the following equalities on
morphisms
\begin{align*}
  (1_{a})^{\dag} &= 1_{a} \\
  f^{\dag\dag} &= f \\
  (f\otimes g)^{\dag} &= f^{\dag}\otimes g^{\dag}
\end{align*}
\end{definition}

Orientation-reversal of cobordisms in the category
$\Catname{HCobord}(n,X)$ satisfies all the axioms for a
$\ast$-structure. In this category, it also follows from the
definitions that if $T$ is a cylinder cobordism then we have the
equality
\begin{equation}
T^\dag = T^{-1}
\end{equation}

Besides the category $\Catname{HCobord}(n,X)$, the category
$\Catname{Hilb}$ of (finite dimensional) Hilbert spaces also has a
$\ast$-structure: It is given by the adjoint of a linear operator.
A symmetric monoidal functor
\begin{equation*}
  \tau:\Catname{HCobord}(n,X)\longrightarrow\Catname{Hilb}
\end{equation*}
satisfying the additional condition $\tau(f^\dag)=(\tau f)^{\dag}$
is called a Unitary homotopy quantum field theory. In this paper
we will not discuss such unitary field theories and instead urge
the reader to go to \cite{JBJD} to see their relevance in physics.

\smallskip
Now let $\Psi:X\longrightarrow Y$ be a continuous map. If $(M,g)$
is an $X$-manifold then $(M,g\Psi)$ is a $Y$-manifold. In the same
way, if $(W,F):(M,g)\longrightarrow(N,h)$ is an $X$-cobordism then
$(W,F\Psi):(M,g\Psi)\longrightarrow(N,h\Psi)$ is a $Y$-cobordism.
It also takes $X$-diffeomorphisms into $Y$-diffeomorphisms in the
obvious manner and a small amount of elementary checking should
convince the reader that these assignments amount to an induced
strict monoidal functor
\begin{equation*}
  \Psi_{\ast}:\Catname{HCobord}(n,X)\longrightarrow\Catname{HCobord}(n,Y)
\end{equation*}

Now let $\Psi$ and $\Phi$ be continuous maps $X\longrightarrow Y$
and $H$ a homotopy between them. If $(M,g)$ is an $X$-manifold
then $(M,g\Psi)$ and $(M,g\Phi)$ are $Y$-manifolds. On the other
hand, the map $g\times 1_{I}H:M\times I\longrightarrow Y$ defined
by
\begin{equation*}
  (m,t)\longmapsto H(g(m),t)
\end{equation*}
is a homotopy between $g\Psi$ and $g\Phi$, and proposition
\ref{thm:modinv} implies that there is a cylinder cobordism
\begin{equation*}
  H_{(M,g)}:\Psi_{\ast}(M,g)\longrightarrow\Phi_{\ast}(M,g)
\end{equation*}

This family of isomorphisms is in fact a natural transformation as
is proved in the next proposition.

\begin{proposition}
The family $H_{(M,g)}$ of cylinder cobordisms is a monoidal
natural isomorphism $\Psi_{\ast}\longrightarrow\Phi_{\ast}$.
\end{proposition}
\begin{proof}
By proposition \ref{thm:isodescription} the naturality question is
reduced to proving the commutativity of the square in figure
\ref{diag:natsquareforhomotopy} for a general $X$-cobordism
$(W,F)$.
\begin{figure} [hbtp]
  $\vcenter{\xymatrix{
    \Psi_{\ast}(M,g) \ar[d]_{\Psi_{\ast}(W,F)} \ar[r]^{H_{(M,g)}} & \Phi_{\ast}(M,g) \ar[d]^{\Phi_{\ast}(W,F)} \\
    \Psi_{\ast}(N,h) \ar[r]^{H_{(N,h)}} & \Phi_{\ast}(N,h)   }
  }$
  \caption{}
  \label{diag:natsquareforhomotopy}
\end{figure}

But since $H_{(M,g)}$ is an isomorphism, commutativity of diagram
\ref{diag:natsquareforhomotopy} is the same as the equality
\begin{equation*}
  \Phi_{\ast}(W,F)=H_{(M,g)}^{-1}\Psi_{\ast}(W,F)H_{(N,h)}
\end{equation*}

The right term consists of two cylinders glued at the ends of $W$
and the map into $Y$ is defined by $\overline{g\times 1_{I}H}\cdot
F\Psi\cdot h\times 1_{I}H$. The proof will be accomplished if we
can deform this map into $1_{\Phi_{\ast}(M,g)}\cdot F\Phi\cdot
1_{\Psi_{\ast}(N,h)}$ by a homotopy that is fixed on the boundary.

In $W\times I$ we define this homotopy simply by
\begin{equation*}
  (w,s)\longmapsto H(F(w),s)
\end{equation*}

Note that if $w=m\in M$ then this reduces to $H(g(m),s)$ and if
$w=n\in N$ then this reduces to $H(h(n),s)$.

For $(N\times I)\times I$ we define the homotopy by
\begin{equation*}
  ((n,t),s)\longmapsto
    \begin{cases}
    H(h(n),t+s) & \text{if $s\leq 1-t$}, \\
    \Phi(h(n)) & \text{if $s\geq 1-t$}.
  \end{cases}
\end{equation*}

This is a continuous function because for $s=1-t$ the first branch
is $H(h(m),1)=\Phi(h(n))$. For $t=0$ this is just $H(h(n),s)$ and
therefore continuity of the whole homotopy is ensured. For $s=1$
this is $\Phi(h(n))$ and therefore the homotopy remains fixed in
this component of the boundary.

We leave to the reader the easy task of supplying the definition
of the homotopy for $(M\times I)\times I$ as well as making the
necessary checks.

As for being a monoidal transformation, since the induced functors
are strict monoidal, this is equivalent to ask that
\begin{equation*}
  H_{(M,g)\amalg(N,h)}=H_{(M,g)}\amalg H_{(N,h)}
\end{equation*}

But this equality follows from the diffeomorphism $I\times(M\amalg
N)\cong(I\times M)\amalg(I\times N)$ and relation
\ref{reln:xdiffrel}.
\end{proof}

We could use the above results to build a $2$-functor between
appropriate $2$-categories, but instead of doing that we turn our
attentions to the fundamental result about the induced functor.

\begin{theorem} \label{thm:weakhomotopyinvariance}
If $\Psi$ is an $(n+1)$-weak homotopy equivalence then the induced
functor $\Psi_{\ast}$ is an equivalence of categories.
\end{theorem}
\begin{proof}
Since our manifolds are differentiable, they admit (at least) one
structure of CW-complexes with cells up to the dimension of the
manifold. Standard obstruction theory says that the map $\Psi$
induces the isomorphisms
\begin{equation*}
\begin{split}
  [M,X] &\cong[M,Y] \\
  [W,X]_{rel\partial W} &\cong[W,Y]_{rel\partial W}
\end{split}
\end{equation*}

It is clear that $\psi$ is an $X$-diffeomorphism iff
$\Psi_{\ast}\psi$ is an $Y$-diffeomorphism, and this fact,
together with an application of the second isomorphism for
$X$-cobordisms, implies that $\Psi_{\ast}$ is full and faithful.
The first isomorphism coupled with an application of proposition
\ref{thm:modinv} implies that $\Psi_{\ast}$ is isomorphism-dense.
\end{proof}

We invite the reader to check, armed with this theorem, that if
the background space is an $(n-1)$-homotopy type we recover
Turaev's axioms, in particular the one for the identity cobordism.

As a corollary of theorem \ref{thm:weakhomotopyinvariance} we have

\begin{corollary}
For a simply-connected space $X$, there is an equivalence of
categories
\begin{equation*}
  \Catname{HCobord}(1,X)\simeq\Catname{HCobord}(1,K(G,2))
\end{equation*}
where $G$ is the abelian group $\pi_{2}X$.
\end{corollary}

In the next section we will give a complete algebraic description
of the category $\Catname{HCobord}(1,K(G,2))$.

Another special case merits our attention: when $X$ is the
one-point space $\ast$. In this situation all the characteristic
maps are redundant, whether in objects or in morphisms, so that
the category $\Catname{HCobord}(n,\ast)$ is a realization of the
category of cobordisms $\Catname{Cobord}(n)$. The terminal map
$X\longrightarrow\ast$ induces a functor
\begin{equation*}
  \Catname{HCobord}(n,X)\longrightarrow\Catname{Cobord}(n)
\end{equation*}

In particular any TQFT can be trivially extended to an HQFT with
target $X$. The other way around is also true. To any HQFT there
is associated an underlying TQFT which is given by restricting to
$X$-manifolds and $X$-cobordisms whose characteristic maps send
everything to the base point $\ast$. In other words there is an
inclusion functor
\begin{equation*}
  \Catname{Cobord}(n)\longrightarrow\Catname{HCobord}(n,X)
\end{equation*}
which is a left-inverse of the above functor.

\smallskip
There is a nice dual aspect here. On the one hand a monoidal
functor on $\Catname{HCobord}(n,X)$ gives us invariants of the
manifolds mapped into $X$. On the other hand we can regard
$\Catname{HCobord}(n,X)$ as a homotopy invariant of $X$ by theorem
\ref{thm:weakhomotopyinvariance} -- we can probe the homotopy type
of $X$ by looking at the manifolds that are mapped into it up to
homotopy.

This dependence on the background space can be made even more
direct by noting the following. Consider the sphere $S^{n}$ of
dimension $n$ with a base-point maintained fixed throughout. The
nth based loop space $\Omega^{n}(X)$ of $X$ is the space
\begin{equation*}
  \Catname{Map}(S^{n},X)
\end{equation*}
of based maps $S^{n}\longrightarrow X$. To each point
$g\in\Omega^{n}(X)$ we can associate the $X$-manifold $(S^{n},g)$,
an object of $\Catname{HCobord}(n,X)$. To a path
$H:I\longrightarrow\Omega^{n}(X)$ such that
\begin{equation*}
  H(0)=g \quad\text{and}\quad H(1)=h
\end{equation*}
we can associate the cylinder cobordism $(S^{n}\times
I,\widetilde{H})$ where $\widetilde{H}$ is a based homotopy, the
image of $H$ by the natural isomorphism
\begin{equation*}
  \Catname{Map}(I,\Omega^{n}(X))\cong\Catname{Map}\left(I,\Catname{Map}(S^{n},X)\right)\cong\Catname{Map}(S^{n}\times I,X)
\end{equation*}

If the paths $H$ and $G$ are homotopic relative to the end-points
then $\widetilde{H}$ and $\widetilde{G}$ are homotopic relative to
the boundary and vice-versa, and this means that, if we denote by
$\pi_{1}(X)$ the fundamental groupoid of $X$, we have

\begin{proposition}
There is an inclusion functor (faithful and injective on objects)
$\pi_{1}\Omega^{n}(X)\longrightarrow\Catname{HCobord}(n,X)$
\end{proposition}

This implies that an HQFT induces a functor
$\pi_{1}\Omega^{n}(X)\longrightarrow\Catname{Mod}(R)$. A functor
such as this is the same thing as a locally constant sheaf on
$\Omega^{n}(X)$ and if we put $R=\Complex$ then it is also the
same thing as a complex flat vector bundle with connection on
$\Omega^{n}(X)$. The invariant associated to a cylinder cobordism
$(S^{n}\times I,H)$ is precisely the parallel transport along the
path $H$ in the loop space $\Omega^{n}(X)$.

It is well known, at least when $X$ is $n$-connected, that a line
bundle in $\Omega^{n}(X)$ is the same thing as an abelian
$n$-gerbe in the base space, so, taking some liberty with the
terminology, we can say that an HQFT in dimension $n+1$ and with
target space $X$ determines a flat $n$-gerbe with connection in
$X$. If $X$ is not $n$-connected then not all $n$-gerbes come from
bundles in the loop space $\Omega^{n}(X)$. We will have more to
say about this case in a forthcoming paper.

\section{HQFT's in dimension $2$ with target space $K(G,2)$.} \label{section:hqftindim2}

In \cite{BT} the authors have considered HQFT's with target space
a simply connected space (or $K(G,2)$ by our results), but here we
go one step further and provide a complete description of
$\Catname{HCobord}(1,X)$ in terms of a universal property. This
description will have as a corollary that the category of HQFT's
with target space $K(G,2)$ is equivalent to the category of
Frobenius algebras with an action of $G$.

The idea in this section is to find a skeleton of
$\Catname{HCobord}(1,X)$ and then characterize it up to
equivalence by a universal property.

\smallskip
We start with some homotopical calculations of a general nature.
Take $X$ to be the Eilenberg-Maclane space $K(G,n+1)$ for some
abelian group $G$ and $n>1$. In this case the objects in
$\Catname{HCobord}(n,K(G,n+1))$ are $n$-dimensional manifolds and
we have
\begin{equation} \label{eqn:iso1}
  [M,K(G,n+1)]\cong H^{n+1}(M;G)\cong 0
\end{equation}

Invoking proposition \ref{thm:modinv} we can, right from the
start, assume that the characteristic maps of our object
$X$-manifolds are trivial, that is, they send everything to the
base point of $X$. But this implies that for an $X$-cobordism $W$
we have
\begin{equation*}
  [W,K(G,n+1)]_{rel\partial W}\cong[(W,\partial
W),(K(G,n+1),\ast)]
\end{equation*}

The set $[W,K(G,n+1)]_{rel\partial W}$ is computed by the
following proposition
\begin{proposition} \label{prop:homotopycalc}
We have a natural isomorphism
\begin{equation*}
  [(W,\partial W),(K(G,n+1),\ast)]\cong H^{n+1}(W,\partial W;G)
\end{equation*}
\end{proposition}
\begin{proof}
The first step is to use the collapsing construction and its
universal property depicted by the diagram in figure
\ref{diag:collapsingconstruct}
\begin{figure}[hbtp]
  $\vcenter{\xymatrix{
  (W,A) \ar[dr] \ar[r] & (W/A,\ast) \ar@{-->}[d]  \\
                & (Y,\ast)}}$
  \caption{}
  \label{diag:collapsingconstruct}
\end{figure}\\
where the map of pairs $(W,A)\longrightarrow(W/A,\ast)$ is the
collapsing map.

This universal property just says that any map of pairs
$(W,A)\longrightarrow(Y,\ast)$ factorizes through a unique map
(drawn in dashed lines) $(W/A,\ast)\longrightarrow(Y,\ast)$. In
terms of isomorphisms it entails
\begin{equation*}
  [(W,A),(K(G,n+1),\ast)]\cong[(W/A,\ast),(K(G,n+1),\ast)]
\end{equation*}
but as $\pi_{1}K(G,n+1)\cong 0$ we have
\begin{equation*}
  [(W/A,\ast),(K(G,n+1),\ast)]\cong[W/A,K(G,n+1)]\cong H^{n+1}(W/A;G)
\end{equation*}

Now use the exact sequence of a pair to get, for $n>1$,
\begin{equation*}
  H^{n+1}(W/A;G)\cong H^{n+1}(W/A,\ast;G)
\end{equation*}

Finally, with $A=\partial W$ and the fact that the inclusion
$\partial W\longrightarrow W$, being a closed cofibration, implies
the natural isomorphism
\begin{equation*}
  H^{n+1}(W/\partial W,\ast;G)\cong H^{n+1}(W,\partial W;G)
\end{equation*}
we get the promised isomorphism of the proposition.
\end{proof}

$W$ is a manifold of dimension $n+1$, so by applying Poincar\'{e}
duality we get
\begin{equation} \label{eqn:iso2}
  H^{n+1}(W,\partial W;G)\cong H_{0}(W;G)
\end{equation}
and by the universal coefficient theorem this amounts to
$H_{0}(W)\otimes G$.

\smallskip
What is a $G$-homology class in degree $0$? It is the simplest
thing one can imagine -- a finite set of points labelled with
elements of $G$. Let us suppose that $W$ is path-connected and
denote by $(w,g)$ a homology class consisting of a point $w$ and
an element $g$. Let $(w',g')$ be another such homology class, then
there is a path $\gamma$ connecting them and we have
\begin{equation*}
  \partial\gamma=g-g'
\end{equation*}

This means that any two such homology classes are the same up to
an element $g$ of $G$. And this implies that
\begin{equation*}
  (w,g)+(w',g')=(w,g+g')=(w',g+g')
\end{equation*}
so that we can denote a homology class by just $g$. This also
implies that we have the

\begin{theorem}[\textbf{Gluing Theorem}] \label{thm:gluethm}
If $(W,g)$ and $(V,g')$ are two composable $X$-cobordisms such
that $W\circ V$ is connected, then we have in
$\Catname{HCobord}(n,K(G,n+1))$ the equality
\begin{equation*}
  (W,g)\circ(V,g')=(W\circ V,g+g')
\end{equation*}
\end{theorem}

\smallskip
Now we apply the above results for $n=1$ to cut down
$\Catname{HCobord}(1,K(G,2))$ to a skeleton and then use the well
known classification theorems for $2$-manifolds to get an
algebraic characterization.

\smallskip
First, notice that we can fix right from the outset a base point
in every connected component of a manifold, since the group of
diffeomorphisms acts transitively in a connected manifold. In
particular, whenever we talk of the circle $S^{1}$ it is always
implicit that it has a base point fixed once and for all.

Next, notice that the isomorphism \ref{eqn:iso1} together with
proposition \ref{thm:modinv} shows that every object in
$\Catname{HCobord}(1,K(G,2))$ is isomorphic to $(M,\ast)$ where
$M$ is a $1$-manifold and $\ast$ is the constant map sending
everything to $\ast\in X$. In particular $K(G,2)$-diffeomorphisms
are nothing but (orientation-preserving) diffeomorphisms. On the
other hand, every $1$-manifold is diffeomorphic to a disjoint
union of copies of the circle. This means that the skeleton is
made up of disjoint unions of $S_{+}^{1}$ and $S_{-}^{1}$, where
$S_{+}^{1}$ ($S_{-}^{1}$) is the circle $S^{1}$ with the standard
positive (negative) orientation.

\smallskip
A diffeomorphism
$\coprod_{i=1}^{n}S^{1}\longrightarrow\coprod_{i=1}^{n}S^{1}$ can
always be written as a permutation composed with
$\coprod_{i}\psi_{i}$ where $\psi_{i}$ are diffeomorphisms
$S^{1}\longrightarrow S^{1}$. Taking into account that the group
of orientation-preserving diffeomorphisms that fix the base point
is contractible, proposition \ref{thm:isotinv} implies that
$\psi_{i}=1_{S^{1}}$.

There is still one diffeomorphism $S_{+}^{1}\longrightarrow
S_{-}^{1}$ which is $-1_{S^{1}}$. By applying proposition
\ref{thm:isodescription} we see that it is equal to a cylinder
cobordism $S_{+}^{1}\longrightarrow S_{-}^{1}$.

\smallskip
Now we make the following definition
\begin{definition}
Let $\mathcal{A}$ be a symmetric monoidal category, where the
monoidal structure is denoted by $\otimes$ and the unit by
$\Field$. A Frobenius object in $\mathcal{A}$ is a tuple
$(a,\mu,\eta,\varrho)$ where $a$ is an object, $\mu$ is a map
$a\otimes a\longrightarrow a$, $\eta$ is a map
$\Field\longrightarrow a$ such that $(a,\mu,\eta)$ is a monoid and
$\varrho$ is a symmetric map $a\otimes a\longrightarrow\Field$
that satisfies the equation depicted in figure
\ref{diag:bilinform}.
\begin{figure}[hbtp]
  $\vcenter{\xymatrix{
    a\otimes a\otimes a \ar[d]_{\mu\otimes 1_{a}} \ar[r]^{1_{a}\otimes\mu} & a\otimes a \ar[d]^{\varrho} \\
    a\otimes a \ar[r]^{\varrho} & \Field   }}$
  \caption{}
  \label{diag:bilinform}
\end{figure}\\
and a non-degeneracy condition in the sense that the adjoint
$a\longrightarrow a^{\ast}$ is an isomorphism.
\end{definition}

\begin{nnexample}
Frobenius objects in the category $\Catname{Mod}(R)$ take the
familiar guise of Frobenius algebras.
\end{nnexample}

\begin{remark}
A Frobenius object is a monoid, so it has a unit. All monoids
(associative algebras, etc.) are assumed to have units. The
respective morphisms always preserve them.
\end{remark}

\begin{notation}
Denote a cobordism by $W_{\epsilon\mu\cdots}$ where
$\epsilon,\mu\in\{-1,+1\}$ denotes the orientation on the
respective circle component of the boundary. If $W$ comes equipped
with a non-trivial homology class $g\in G$ we just write
$W_{\epsilon\mu\cdots}(g)$.

Denote by $P_{--+}$ the pants viewed as a cobordism
$S_{+}^{1}\amalg S_{+}^{1}\longrightarrow S_{+}^{1}$, by $B$ (or
$C_{--}$) the cylinder $I\times S^{1}$ viewed as a cobordism
$S_{+}^{1}\amalg S_{+}^{1}\longrightarrow\emptyset$, by $D_{+}$
the disk as a cobordism $\emptyset\longrightarrow S_{+}^{1}$, and
by $C_{-+}$ the cylinder viewed as an endomorphism cobordism
$S_{+}^{1}\longrightarrow S_{+}^{1}$. Two other important
cylinders are $C_{+-}:\emptyset\longrightarrow S_{+}^{1}\amalg
S_{-}^{1}$ which is just the map $\eta_{S_{+}^{1}}$, and
$C_{+-}:S_{-}^{1}\amalg S_{+}^{1}\longrightarrow\emptyset$ which
is the map $\varepsilon_{S_{+}^{1}}$. These morphisms are the
basic morphisms that generate the whole category
$\Catname{HCobord}(1,K(G,2))$.
\end{notation}

\smallskip
It is well-known that the pants
\begin{equation*}
  P_{--+}:S_{+}^{1}\amalg S_{+}^{1}\longrightarrow S_{+}^{1}
\end{equation*}
the disk $D_{+}$ and the cylinder $B$ make $S_{+}^{1}$ into a
Frobenius object (see \cite{LA} for example). The new structure is
related to the cylinders $C_{-+}(g)$ which are endomorphisms
$S_{+}^{1}\longrightarrow S_{+}^{1}$. An application of theorem
\ref{thm:gluethm} shows that the assignment $g\longmapsto
C_{-+}(g)$ gives a homomorphism
\begin{equation*}
  G\longrightarrow End(S_{+}^{1})
\end{equation*}
and the gluing theorem yields the commutative
diagrams\footnote{The reader is urged to draw the corresponding
pictures and see the correctness of the assertions.} in figure
\ref{diag:frobact}.
\begin{figure}[hbtp]
  $\vcenter{\xymatrix{
    S_{+}^{1}\amalg S_{+}^{1} \ar[dd]_{C_{-+}(g)\amalg C_{-+}} \ar[rr]^{C_{-+}\amalg C_{-+}(g)} \ar[dr]^{P_{--+}} & & S_{+}^{1}\amalg S_{+}^{1} \ar[dd]^{P_{--+}} \\
    & S_{+}^{1} \ar[dr]_{C_{-+}(g)} & \\
    S_{+}^{1}\amalg S_{+}^{1} \ar[rr]^{P_{--+}} & & S_{+}^{1}   }}$
  $\vcenter{\xymatrix@C+30pt{
    S_{+}^{1}\amalg S_{+}^{1} \ar[d]_{C_{-+}(g)\amalg C_{-+}} \ar[r]^{C_{-+}\amalg C_{-+}(g)} & S_{+}^{1}\amalg S_{+}^{1} \ar[d]^{B} \\
    S_{+}^{1}\amalg S_{+}^{1} \ar[r]^{B} & \emptyset   }}$
  \caption{}
  \label{diag:frobact}
\end{figure}

\begin{definition}
A $G$-Frobenius object in a monoidal category $\mathcal{A}$ is a
Frobenius object $a$ endowed with an action $G\longrightarrow
End(a)$ satisfying the axioms given by the commutative diagrams in
figure \ref{diag:frobact} with $S_{+}^{1}$ replaced with $a$,
etc..
\end{definition}

In the next example we translate this abstract definition into
more concrete terms by looking at the monoidal category
$\Catname{Mod}(R)$.

\begin{nnexample}
$G$-Frobenius objects on the category of $R$-modules
$\Catname{Mod}(R)$ will be called $G$-Frobenius algebras. If we
denote the product in a $G$-Frobenius algebra $A$ by
concatenation, the bilinear form by $(,)$ and the action of an
element $g\in G$ by $\cdot$ then the commutative diagrams in
figure \ref{diag:frobact} give the equations
\begin{equation*}
  v(g\cdot w)=g\cdot(vw)=(g\cdot v)w
\end{equation*}
and
\begin{equation*}
  (v,g\cdot w)=(g\cdot v,w)
\end{equation*}

The fact that $A$ is a unital algebra and the above two equations
entail some very strong conditions on the action of $G$. Because
$A$ is unital we have
\begin{equation*}
  g\cdot v=g\cdot(1v)=(g\cdot1)v
\end{equation*}
which means that the action factorizes through the map
$G\longrightarrow A$ given by $g\longmapsto g\cdot 1$.

This map is a homomorphism since
\begin{equation*}
  (gh)\cdot 1=g\cdot(h\cdot 1)=g\cdot(1(h\cdot 1))=(g\cdot 1)(h\cdot 1)
\end{equation*}

Also, $g\cdot 1$ is in the center of $A$ because
\begin{equation*}
  (g\cdot 1)v=g\cdot(1v)=g\cdot(v1)=v(g\cdot 1)
\end{equation*}

Below we will show that a $(1+1)$-HQFT with target space $K(G,2)$
is essentially determined by the $G$-Frobenius algebra on the
circle. The above computations give the important fact that the
action factorizes through a map of $G$ into the center of $A$.
\end{nnexample}

\smallskip
$S_{+}^{1}$ is a $G$-Frobenius object in
$\Catname{HCobord}(1,K(G,2))$. Below we outline a demonstration of
the fact that is is the universal such object, that is, given a
$G$-Frobenius object $a$ in a symmetric monoidal category
$\mathcal{A}$ there is a unique symmetric monoidal functor
$\Catname{HCobord}(1,X)\longrightarrow\mathcal{A}$ taking
$S_{+}^{1}$ to $a$ and $S_{-}^{1}$ to $a^{\ast}$. In other words,
we have the following elegant theorem

\begin{theorem}
The monoidal category $\Catname{HCobord}(1,K(G,2))$ is equivalent
to the free symmetric monoidal category with strict duals on one
$G$-Frobenius object.
\end{theorem}
\begin{proof}
Part of the proof was done above when we have showed that any
object is isomorphic to a disjoint union of copies of $S_{-}^{1}$
and $S_{+}^{1}$ and shown that every $X$-diffeomorphism is either
a symmetry or can be turned into a cylinder cobordism.

Note that by the classification theorems of $2$-manifolds, plus
the gluing theorem \ref{thm:gluethm} and the above commutative
diagrams in figure \ref{diag:frobact}, for any $2$-manifold that
has a non-trivial homology class $g$ in it, we can always push the
$g$ class to a cylinder.

Also by the classification theorems of $2$-manifolds, every such
manifold with a given orientation can be decomposed into a gluing
and disjoint union of the basic morphisms listed above. On the
other hand every such manifold, viewed as a cobordism in some
particular way, can be decomposed in terms of these basic
morphisms and their duals, that is, compositions of the basic
morphisms with the pairings $\eta$ and $\varepsilon$.

Of course, for a given cobordism there may be more than one way of
writing it as a composition of the basic cobordisms so we have to
ensure that all the relations are covered by the algebraic
relations of a Frobenius object and the relations of a symmetric
monoidal category. But this is a purely topological fact
pertaining to Cerf theory and that the reader can see in
\cite{LA}.
\end{proof}

There are two immediate corollaries worth mentioning.

\begin{corollary}
The category $\Catname{HQFT}(1,K(G,2))$ is equivalent to the
category of $G$-Frobenius algebras.
\end{corollary}

If we take $G$ to be the trivial group $0$, $K(0,2)$ is just the
one-point space and we have

\begin{corollary}
The category $\Catname{Cobord}(1)$, the category of
$(1+1)$-cobordisms, is equivalent to the free symmetric monoidal
category with strict duals on one Frobenius object.
\end{corollary}

\section{State-sum models for $2d$-HQFT's}

We want to modify the state sum models for a TQFT by incorporating
homotopical data. To do that we will in the first place review the
results in \cite{FHK} since they are a necessary springboard for
our computations.

\smallskip
Let $(M,T)$ be a closed triangulated surface. We will work with
the dual graph $T^{\ast}$ of the triangulation $T$, and instead of
explaining what that is we present a picture of the dual graph of
a triangle in figure \ref{grph:dualgraph}.
\begin{figure}[hbtp]
  $\vcenter{\xymatrix{
    \ar@{--}[ddr] & & \ar@{--}[ddl] \\
     & *[o][F]{} \ar@{-}[ddl] \ar@{-}[ddr] & \\
     & *+[o][F]{} \ar@{--}[dd] & \\
    *[o][F]{} \ar@{-}[rr] & & *[o][F]{} \\
     & &   }}$
  \caption{}
  \label{grph:dualgraph}
\end{figure}

The vertices of the dual graph are displayed as big round empty
circles, that will be filled later with elements of an abelian
group $G$. The edges of the dual graph are drawn in dashed lines.

\begin{definition}
A flag is a pair $(e,v)$, where $e$ is an edge of the dual graph
and $v$ is one of the vertices to which to edge is incident.

A coloring of $T$ is an assignment of an element $i\in
\{0,\cdots,n\}$ to each flag of the dual graph.
\end{definition}

The reader should note that each internal edge is colored twice,
since it appears in two flags. In terms of the original
triangulation, a coloring corresponds to first ``exploding'' the
triangulation, so that each edge is duplicated, and then make an
assignment of $i\in \{0,\cdots,n\}$ to each edge.

A coloring can be seen as filling our surface with fields. The
dual graph can be interpreted as a kind of Feynman graph, and what
we will do, is apply a Feynman-like rule to compute its amplitude.

\textbf{Rule:} To each vertex with incident edges colored with
$i$,$j$ and $k$ we assign a (complex) number $C_{ijk}$, the vertex
contribution, and to each internal line colored with $k$ and $k'$
we assign a propagator $g^{kk'}$. Then we multiply all these
numbers and sum over all possible colorings, which amounts to
contracting all the indices, to obtain the partition function
$Z(M,T)$.

Alas, there is still one little detail to take care of: what is
the order in which the indices appear in $C_{ijk}$? To answer
that, note that an orientation of $M$ induces an orientation on
each $2$-simplex in $T_{(2)}$, which in its turn induces an
orientation on its boundary. This induced orientation gives a
cyclic order to the edges incident to a vertex in the dual graph
and this is the order we use to write down $C_{ijk}$. Since the
order is only cyclic we impose that the $C_{ijk}$ remains
invariant when we apply a cyclic permutation to its indices. The
same reasoning prompts us to impose that $g^{kk'}$ is symmetric in
the indices.

In figure \ref{grph:graphcolor} we have a picture of such a
situation. Here the order on the edges is to be identified with
the alphabetic order of the colorings.
\begin{figure}[hbtp]
  $\vcenter{\xymatrix{
    \ar@{--}[dr]^{j} & & *[o][F]{} \ar@{-}[dll] \ar@{-}[dd] \ar@{-}[drr] & & \ar@{--}[dl]_{l} \\
    *[o][F]{} \ar@{-}[drr] & *+[o][F]{} \ar@{--}[rr]_>>>>{k'}^<<<<{k} \ar@{--}[dl]^{i} & & *+[o][F]{} \ar@{--}[dr]_{m} & *[o][F]{} \ar@{-}[dll] \\
     & & *[o][F]{} & &   }}$
  \caption{}
  \label{grph:graphcolor}
\end{figure}
To this colored graph corresponds the expression
\begin{equation*}
  C_{ijk}g^{kk'}C_{k'lm}=C_{ij}^{\phantom{ij}k}C_{klm}
\end{equation*}

The authors of \cite{FHK} go on to show that taking the free
linear space $A$ generated by $\{0,\ldots,n\}$, and imposing the
constraint that $Z(M,T)$ remains invariant under (a suitable
version of) the Pachner moves that express the topological
invariance of the triangulation, $C_{ij}^{\phantom{ij}k}$ are the
coefficients for a bilinear associative multiplication and
$g^{ik}$ are the coefficients of a non-degenerate bilinear form.
In particular Topological state-sums (or Topological lattice field
theories) are in bijection with semi-simple associative algebras
$A$.

Now we want to generalize this state sum in order to include
homotopical information on $M$, that is, we want invariants
$Z(M,T,g)$ where $g$ is now a homotopy class of maps
$M\longrightarrow X$. By theorem \ref{thm:weakhomotopyinvariance}
only the nth weak homotopy type of $X$ matters, so choose a
simplicial approximation $\widetilde{X}$ of $X$, that is, a
simplicial complex $\widetilde{X}$ with a map
$\widetilde{X}\longrightarrow X$ inducing isomorphisms on all
homotopy groups\footnote{Here $\widetilde{X}$ denotes the
geometric realization of the complex. We are making no notational
distinction between a simplicial complex and its geometric
realization.}. Now by the simplicial approximation theorem any map
$M\longrightarrow\widetilde{X}$ is homotopic to a simplicial map
and since $X$ is path connected we can choose $\widetilde{X}$ to
have only one vertex $\ast$ so that every vertex of $M$ gets sent
by $g$ into $\ast$.

This gives the combinatorial picture. Our triangulation is right
from the start ``colored'' with cells of the simplicial complex
$\widetilde{X}$ which, in line with the remarks made in the
introduction and elsewhere, we interpret as matter fields.

To get an invariant we color the flags of the dual graph of the
triangulation $T$ with, this in dimension $1+1$, elements from the
finite set $\{0,\ldots,n\}$, which, once again taking on the hints
from the introduction, can be seen as geometries. The rule for the
state sum is the same except that we make the following
replacements in the vertex and edge contribution
\begin{equation*}
\begin{split}
  C_{ijk} &\longmapsto C_{ijk}(g) \\
  g^{kk'} &\longmapsto g^{kk'}(\gamma)
\end{split}
\end{equation*}
where $\gamma$ ($g$) is a $1$-cell ($2$-cell) of $\widetilde{X}$.

The proof of the invariance now goes in two steps. The first is to
prove that if we change $g$ to a homotopic map $g'$ the state sum
remains invariant and the second step is to prove that the sum is
invariant by changing the triangulation $T$ to another equivalent
triangulation $T'$. Let us suppose that the first step has been
achieved then equivalent triangulations $T$ and $T'$ are related
by a finite number of Pachner moves. On the other hand the region
$P$ where the Pachner move is performed is a subcomplex and
therefore the inclusion $P\longrightarrow T$ is a cofibration. $P$
is also contractible and therefore the quotient map
$M\longrightarrow M/P$ is a homotopy equivalence which implies
that in the combinatorial description of the map $g$ we can assume
that it is trivial in $P$. Rerunning the proof in \cite{FHK} we
realize that $C_{ij}^{\phantom{ij}k}(\ast)=C_{ij}^{\phantom{ij}k}$
are the coefficients of an associative multiplication in the free
linear space generated by $\{0,\ldots,n\}$ and
$g^{kk'}(\ast)=g^{kk'}$ are the coefficients of a bilinear
non-degenerate metric in this same space. In particular,
topological invariance is achieved.

So now what we are left to do is to prove homotopical invariance,
or better said, to discover what constraints homotopical
invariance puts on the factors $C_{ijk}(g)$, $g^{kk'}(\gamma)$.
Since our target space is $K(G,2)$ we can choose the simplicial
approximation to have only $1$-cell and this means that the
propagators $g^{kk'}$ are automatically the coefficients of the
metric.

We have the isomorphism
\begin{equation*}
  [M,K(G,2)]\cong H^{2}(M;G)
\end{equation*}

From the universal coefficient theorem we know that
\begin{equation*}
  H^2(M;G)\cong Hom(H_{2}(M);G)
\end{equation*}
and since $M$ is orientable, by fixing generators this implies
that $H^2(M;G)\cong H_{2}(M;G)\cong G$ so we can view $g$ as an
assignment of an element of $G$ to each $2$-simplex of the
triangulation. The last isomorphism just says that we can shift
the $g\in G$ from a $2$-simplex to another $2$-simplex colored
with $h$ if we sum them. In terms of dual graphs this leads to the
two constraints that we present in a pictorial form in figure
\ref{graph:basicmove}.
\begin{figure}[hbtp]
  $\vcenter{\xymatrix{
     \ar@{-}[dr] & & & \ar@{-}[dl] & & \ar@{-}[dr] & & & \ar@{-}[dl] \\
    & *++[o][F]{g} \ar@{-}[dl] \ar@{-}[r] & *++[o][F]{h} \ar@{-}[dr] & & \rightleftharpoons & & *++[o][F]{\phantom{g}} \ar@{-}[dl] \ar@{-}[r] & *+<3mm,3mm>[o][F]{gh} \ar@{-}[dr] & \\
    & & & & & & & &   }}$
  $\vcenter{\xymatrix{
     \ar@{-}[dr] & & & \ar@{-}[dl] & & \ar@{-}[dr] & & & \ar@{-}[dl] \\
    & *++[o][F]{g} \ar@{-}[dl] \ar@{-}[r] & *++[o][F]{h} \ar@{-}[dr] & & \rightleftharpoons & & *+<3mm,3mm>[o][F]{gh} \ar@{-}[dl] \ar@{-}[r] & *++[o][F]{\phantom{g}} \ar@{-}[dr] & \\
    & & & & & & & &   }}$
  \caption{}
  \label{graph:basicmove}
\end{figure}

The empty circles mean that the vertex is filled with the unit
element $1\in G$. Notice also how we have used multiplicative
notation.

As in \cite{FHK} we will derive the constraints that $C_{ijk}(g)$
has to satisfy under these homotopical moves. The constraints of
figure \ref{graph:basicmove} on $C_{ijk}(g)$ written in equations
mean
\begin{equation} \label{eqn:basicmove}
  C_{ij}^{\phantom{ij}k}(g)C_{klm}(h)=C_{ij}^{\phantom{ij}k}C_{klm}(gh)=C_{ij}^{\phantom{ij}k}(gh)C_{klm}
\end{equation}

For the same reasons as explained above, we can suppose that
$C_{ijk}(g)$ is cyclic in the indices.

For $g=1$ we have the useful equation
\begin{equation} \label{eqn:singlemove}
  C_{ij}^{\phantom{ij}k}C_{klm}(h)=C_{ij}^{\phantom{ij}k}(h)C_{klm}
\end{equation}

We can always assume that the element $0$ is the multiplicative
identity of the algebra $A$; if not we just make a suitable change
of basis. Now, the game will consist in starting from the above
equations, putting one of the free indices to $0$ and seeing where
we can get. We start with $i=0$ in the equations
\eqref{eqn:basicmove} and \eqref{eqn:singlemove} to get
\begin{equation} \label{eqn:eqn1}
  C_{0j}^{\phantom{0,j}k}(g)C_{klm}(h)=C_{jlm}(gh)
\end{equation}
\begin{equation} \label{eqn:eqn2}
  C_{0j}^{\phantom{0,j}k}(g)C_{klm}=C_{jlm}(g)
\end{equation}

From equation \eqref{eqn:eqn1}, by raising the free index $m$ and
with $l=0$ we get
\begin{equation} \label{eqn:compeqn}
  C_{j0}^{\phantom{j0}m}(gh)=C_{0j}^{\phantom{0j}k}(g)C_{k0}^{\phantom{k0}m}(h)
\end{equation}

Putting $h=1$ we get
\begin{equation} \label{eqn:symlow}
C_{j0}^{\phantom{j0}m}(g)=C_{0j}^{\phantom{0j}k}(g)C_{k0}^{\phantom{k0}m}=C_{0j}^{\phantom{0j}k}(g)\delta_{k}^{m}=C_{0j}^{\phantom{0j}m}(g)
\end{equation}

If we view $C_{0j}^{\phantom{0j}k}(g)$ as the matrix elements of
an action of $G$ in our initial algebra, employing the notation
introduced in section \ref{section:hqftindim2} we have
\begin{equation*}
  C_{j0}^{\phantom{j0}m}(g)=(g\cdot j)_{m}
\end{equation*}
where the subscript $m$ means the component $m$ in the chosen
basis $\{0,\ldots,n\}$. Now, equation \eqref{eqn:compeqn} together
with \eqref{eqn:symlow} states the homomorphism property. Equation
\eqref{eqn:eqn2} can be written as
\begin{equation*}
  C_{jl}^{\phantom{jl}m}(g)=((g\cdot j)l)_{m}
\end{equation*}

Back to equation \eqref{eqn:eqn1} with $h=1$ and $j=0$ we get
\begin{equation} \label{eqn:actionfact}
  C_{0l}^{\phantom{0l}m}(g)=C_{00}^{\phantom{00}k}(g)C_{kl}^{\phantom{kl}m}
\end{equation}
and this equation, in the more invariant notation, is
\begin{equation*}
  (g\cdot l)_{m}=((g\cdot 1)l)_{m}
\end{equation*}

Notice that using cyclicity in \eqref{eqn:actionfact} we would get
\begin{equation*}
  C_{0m}^{\phantom{0m}l}(g)=C_{00}^{\phantom{00}k}(g)C_{mk}^{\phantom{mk}l}
\end{equation*}
which means
\begin{equation*}
  (g\cdot m)_{l}=(m(g\cdot 1))_{l}
\end{equation*}

Now plug these in equation \eqref{eqn:compeqn}, setting $j=0$, to
get
\begin{equation} \label{eqn:homeqn}
  C_{00}^{\phantom{00}m}(gh)=C_{00}^{\phantom{00}k}(g)C_{00}^{\phantom{00}l}(h)
C_{lk}^{\phantom{lk}m}
\end{equation}
which means that the map $g\longmapsto g\cdot 1$ is a
homomorphism
\begin{equation*}
  (gh\cdot 1)_{m}=((g\cdot 1)(h\cdot 1))_{m}
\end{equation*}

We can draw several conclusions from these equations. First, the
action is Frobenius since (dropping the subscript components)
\begin{equation*}
  (g\cdot m)l=((g\cdot 1)m)l=(g\cdot 1)(ml)=g\cdot(ml)
\end{equation*}

Second, the map $g\longmapsto g\cdot 1$ falls on the center of
$A$, since
\begin{equation*}
  (g\cdot 1)k=g\cdot(1k)=g\cdot(k1)=k(g\cdot 1)
\end{equation*}

The corollary of these computations is the statement of the
theorem that sums up the state of affairs for homotopical state
sums in dimension $1+1$ and target space $K(G,2)$.

\begin{theorem}
Any $(1+1)$-dimensional topological state sum built from a
semi-simple algebra $A$ automatically enriches to an homotopical
state sum (HSS for short) with target space $K(Z(A)^{\ast},2)$,
where $Z(A)^{\ast}$ are the invertible elements of the center
$Z(A)$ of $A$. Reducing the $Z(A)^{\ast}$-HSS  to a group $G$
amounts to finding a homomorphism $\varphi:G\longrightarrow
Z(A)^{\ast}$.

Every HSS can be built in this way, in other words, $K(G,2)$-HSS
are in bijective correspondence with pairs $(A,\varphi)$ where $A$
is a semi-simple algebra and $\varphi$ is a homomorphism
$G\longrightarrow Z(A)^{\ast}$.
\end{theorem}

\section{Conclusions.}

In this paper we have characterized $\Catname{HCobord}(1,K(G,2))$
as the free symmetric monoidal category with strict duals on one
$G$-Frobenius object and have given the corresponding state-sum
models. In a forthcoming paper we will give all the details for
the case of a general homotopy $2$-type. In this case the action
of $G$ is, in rough terms, replaced by the functorial action of
$\pi_{1}\Omega(X)$ so that, in the state sum models, the important
ingredient is a flat $Z(A)^{\ast}$-gerbe in $X$. In other words
topological matter in $1+1$ dimensions is the same thing as a
gerbe in $X$.

It is very interesting that gerbes are also appearing in string
theory and related subjects, see, for example, \cite{JK} and
\cite{YZ}.

In higher dimensions we expect to observe the same
categorification pattern already recognized in TQFT's. It is well
known that in three dimensions topological state sums are built
from certain monoidal categories (see for example \cite{TQ}) and
that in four dimensions this role is played by certain monoidal
$2$-categories (see \cite{MM}). To enrich these state sum models
to homotopical ones we expect that what is needed is a functor
from the higher-dimensional fundamental groupoids into the
(generalized) center of these categories.

Besides going up the dimensional ladder, we wish to mention two
further avenues of research. In the paper \cite{MMRP} the authors
have shown that when $X$ is simply connected, gerbes are
essentially the same thing as morphisms from the thin fundamental
group of the loop space into the structure group of the gerbe.
This raises the question of wether it is possible to define Thin
Homotopy Quantum Field theories, by using thin homotopy classes of
maps.

The second avenue is that as there is a classifying space $BG$ of
$G$-bundles, P. Gajer (\cite{GA}) has introduced a classifying
space for (abelian) gerbes. Starting with these classifying
spaces, the constructions made in \cite{FQ} for TQFT's should be
generalized to HQFT's.

\section*{Appendix: The category $\Catname{HCobord}(n,X)$}
\setcounter{figure}{0}

In this appendix we give the details of the construction of the
category of homotopy cobordisms $\Catname{HCobord}(n,X)$. This
construction works in every dimension and encompasses all the
desired features for such a category including the action of the
diffeomorphisms.

The idea of the construction is in fact pretty straightforward.
Build a category in which the objects are $X$-manifolds and where
the morphisms are not only $X$-cobordisms but also
$X$-diffeomorphisms. The relations for composing them are given by
reading off the axioms for an HQFT stripped of the HQFT $\tau$.
While the idea is simple, nevertheless there are some technical
details to surmount. The main problem, as the reader might have
guessed, is with the composition of cobordisms. It's not
associative and the unit cobordisms are not units at all. So some
quotient is needed to get things right. The overkill procedure of
taking equivalence classes of manifolds does not work because it
would not only kill the action of the diffeomorphisms, it would
also make the simple idea of a map into some background space $X$
an undefinable concept. In fact, one of the key ideas is that the
construction proceeds not by reduction -- taking equivalence
classes of manifolds -- but by enlargement -- by adding
diffeomorphisms as arrows.

We want to stress that this construction is purely algebraic. The
role of topology and geometry ends when we have given a name to
our objects. This can either be good or bad depending on one's own
preferences but it certainly has the advantage of robustness and
portability in that it can be applied to a large number of similar
situations where we have a functorial-like invariant of some
geometro-topological object.

We will divide the construction in a number of elementary steps
starting with

\smallskip
\textbf{Step 1: A preliminary construction.}

Before sinking our teeth in $\Catname{HCobord}(n,X)$ let us make
an excursion and explain a simple quotient construction that we
will need later. The setup is the following: we have a category
$\mathcal{A}$ and a subcategory $\mathcal{G}$ which is also a
groupoid.

Define the following relation on the objects of $\mathcal{A}$

\begin{equation*}
  a\sim b \quad\text{iff}\quad \exists g\in\mathcal{G}:a\longrightarrow b
\end{equation*}

Because $\mathcal{G}$ is a groupoid this is clearly an equivalence
relation. Define also the following relation in the arrows of
$\mathcal{A}$: $f\sim g \quad\text{iff}\quad
\exists\psi,\phi\in\mathcal{G}$ such that the square in figure
\ref{diag:relsquare} is commutative.
\begin{figure}[hbtp]
  $\vcenter{\xymatrix{
    a \ar[d]_{f} \ar[r]^{\psi} & a' \ar[d]^{g} \\
    b \ar[r]^{\phi} & b'   }}$
  \caption{}
  \label{diag:relsquare}
\end{figure}

Once again, this is an equivalence relation because $\mathcal{G}$
is a groupoid. So, define $\mathcal{A}/\!\!\sim$ having as objects
the equivalence classes of objects by the above relation, and as
arrows the equivalence classes of arrows of $\mathcal{A}$.

As for the source and target functions we simply have
\begin{align*}
  s[f] &=[s(f)] \\
  t[f] &=[t(f)]
\end{align*}
and it is easy to see that this definition does not depend on the
representatives chosen.

The operation of composition is the one coming from $\mathcal{A}$
and the question is now under what conditions is this well
defined. The answer is provided by the following

\begin{nnlemma}
If the groupoid $\mathcal{G}$ is thin (it has at most one arrow
between any two objects) then we have a well-defined quotient.

In this case the quotient functor
$\pi:\mathcal{A}\longrightarrow\mathcal{A}/\!\!\sim$ is an
equivalence of categories.
\end{nnlemma}
\begin{proof}
We need to prove that the relation on arrows is in fact a
congruence. So, let $f_1\sim g_1$ and $f_2\sim g_2$, then we have
the commutative squares in figure \ref{diag:compsquare}.
\begin{figure}[hbtp]
  $\vcenter{\xymatrix{
    a \ar[d]_{f_1} \ar[r]^{\psi} & a' \ar[d]^{g_1} \\
    b \ar[r]^{\phi} & b'   }}$
  $\vcenter{\xymatrix{
    b \ar[d]_{f_2} \ar[r]^{\phi'} & b' \ar[d]^{g_2} \\
    c \ar[r]^{\varphi} & c'   }}$
  \caption{}
  \label{diag:compsquare}
\end{figure}

Since $\mathcal{G}$ is thin we have $\phi=\phi'$. So, compose the
two squares in diagram \ref{diag:compsquare} top to bottom to
obtain the desired $f_1 f_2\sim g_1 g_2$.

To see that $\pi$ is an equivalence note that if $f$ and $g$ in
diagram \ref{diag:relsquare} are both in $\mathcal{A}(a,b)$ then
because $\mathcal{G}$ is thin we have $\psi=1_a$ and $\phi=1_b$
but this immediately implies that the map
$\mathcal{A}(a,b)\longrightarrow\mathcal{A}/\!\!\sim(\pi(a),\pi(b))$
is bijective. Finally, the fact that $\pi$ is surjective on
objects implies that it is an equivalence.
\end{proof}

If in diagram \ref{diag:relsquare} we put $\psi$ and $g$ equal to
the identity, and also $f=\psi$ and $\phi=\psi^{-1}$ we get that
$\psi\sim 1$ and therefore the structural isomorphisms go into
identities. The whole idea of the construction is to identify
pairs of objects for which there is a canonical (e.g. in
$\mathcal{G}$) isomorphism connecting them. But to get this right
the arrows of $\mathcal{G}$ must go into identities, and this is
only possible if the only automorphisms present in $\mathcal{G}$
are the trivial ones.

We also state the universal property of this quotient.

\begin{nnlemma}
If $\mathcal{A}$ is a category in the above conditions and
$F:\mathcal{A}\longrightarrow\mathcal{B}$ such that
\begin{equation*}
  a\sim b\Longrightarrow F(a)=F(b)
\end{equation*}
and the same for arrows, that is
\begin{equation*}
  f\sim g\Longrightarrow F(f)=F(g)
\end{equation*}
then there exists a unique functor
$\mathcal{A}/\mathcal{G}\longrightarrow\mathcal{B}$ such that the
diagram in figure \ref{diag:univpropofquocient} is commutative.
\begin{figure}[hbtp]
  $\vcenter{\xymatrix{
  \mathcal{A} \ar[d]_{\pi} \ar[r]^{F} & \mathcal{B}\\
  \mathcal{A}/\mathcal{G} \ar@{-->}[ur] &  }
  }$
  \caption{}
  \label{diag:univpropofquocient}
\end{figure}
\end{nnlemma}
\begin{proof}
Just define the functor, called $F'$, say, as $F'[a]=F[a]$ and
$F'[f]=F[f]$. The above conditions on $F$ guarantee that this is a
sound definition and the rest follows.
\end{proof}

This gives the natural (in $\mathcal{B}$) isomorphism of
\emph{sets}
\begin{equation*}
  \Catname{Cat'}(\mathcal{A},\mathcal{B})\cong\Catname{Cat}(\mathcal{A}/\mathcal{G},\mathcal{B})
\end{equation*}
where $\Catname{Cat}$ denotes the set of functors and
$\Catname{Cat'}$ denotes the set of functors satisfying the
conditions of the lemma. But it is an easy exercise in category
theory that this lifts to an isomorphism of \emph{categories}.

Finally note that since $\pi$ is an equivalence, if $\mathcal{A}$
is a monoidal category we can just push the monoidal structure of
$\mathcal{A}$ to $A/\!\!\sim$ to obtain a monoidal category. In
this case the functor $\pi$ is tautologically a (strict) monoidal
functor.

\smallskip
\textbf{Step 2: Notation and structures.}

Recall that we have the category $\Catname{Diff}(n,X)$ having as
objects $n$ dimensional $X$-manifolds and morphisms
$X$-diffeomorphisms. This category is symmetric monoidal by
disjoint union of $X$-manifolds and the unit is the empty manifold
$\emptyset$. We will denote the associator by $\alpha$ (possibly
with subscripts), the unital structural maps by
$l_{M}:M\amalg\emptyset\longrightarrow M$ and
$r_{M}:\emptyset\amalg M\longrightarrow M$ and the symmetry by
$\sigma_{M,N}:M\amalg N\longrightarrow N\amalg M$.

\begin{notation}
Here and below, whenever it is convenient for the sake of clarity,
we will omit the reference to the characteristic maps.
\end{notation}

Now introduce the category $\Catname{C}(n,1;X)$ whose objects are
$X$-cobordisms bounding $n$ dimensional $X$-manifolds, and the
morphisms are $X$-diffeomorphisms between them. It is, just as
$\Catname{Diff}(n,X)$, a symmetric monoidal category for disjoint
union. Denote by $\Gamma_{W,V,U}$, $L_W$, $R_W$ and $\Sigma_{W,V}$
the structural isomorphisms.

Note that $\Gamma_{W,V,U}$ when restricted to the boundary is
precisely the associator $\alpha$, so that, applying axiom
\ref{axm:ax6} of the definition of an HQFT we get the commutative
diagram in figure \ref{diag:cobordassoc}.
\begin{figure}[hbtp]
  $\vcenter{\xymatrix@C+30pt{
    \tau(M\amalg(N\amalg P)) \ar[d]_{\tau(W\amalg(V\amalg U))} \ar [r]^{\tau(\alpha_{M,N,P})} & \tau((M\amalg N)\amalg P) \ar[d]^{\tau((W\amalg V)\amalg U)} \\
    \tau(M'\amalg(N'\amalg P')) \ar [r]^{\tau(\alpha_{M',N',P'})} & \tau((M'\amalg N')\amalg P')   }}$
  \caption{}
  \label{diag:cobordassoc}
\end{figure}

The same kind of relation exists for the maps $L_W$, $R_W$ and
$\Sigma_{W,V}$.

\smallskip
As for gluing $X$-cobordisms we remark that we have a natural
$X$-diffeomorphism
\begin{equation*}
  \Delta_{W,V,U}:W\circ(V\circ U)\longrightarrow (W\circ V)\circ U
\end{equation*}
which is the identity when restricted to the boundary. The key
property of these isomorphisms is that they satisfy the same
pentagon equation for an associator in a monoidal category.

We note, however, that there are no \emph{natural} isomorphisms
$(I\times M,1_g)\circ V\longrightarrow V$ or $W\circ(I\times
M,1_g)\longrightarrow W$. This is the source of some nuisance that
will force us to work with categories without units -- we call
them simply wu-categories. As a saving grace, wu-categories are
not that much different from plain old categories. The same
terminology (functor, natural transformation, \ldots) and the same
constructions (the free wu-category on a graph, \ldots) apply
equally well.

\smallskip
\textbf{Step 3: Making the gluing associative.}

Consider the category $\Catname{C}(n,1,X)$ whose objects are
$X$-cobordisms and the morphisms are $X$-diffeomorphisms between
them. The basic idea is to consider the groupoid $\mathcal{G}$ of
$\Catname{C}(n,1,X)$ generated by the associators
$\Delta_{W,V,U}$, use Maclane's coherence theorem to show that it
is thin and then apply the above construction-lemma. Unfortunately
the coherence theorem cannot be applied as is because, and we
quote \cite{MAC}, \textit{``two apparently or formally different
vertices of such a diagram might become equal in a particular
monoidal category, in such a way as to spoil commutativity.''} So
we have to resort to some formal trick and then, in the end, make
sure that it makes no real difference.

The way it was done in \cite{MAC} was to consider words in
objects, or lists of objects with chosen parenthesations. We
choose something similar but more suitable to our purposes:
redefine the category $\Catname{C}(n,1,X)$ as having as objects
pairs $(W,l)$ where $W$ is an $X$-cobordism and $l$ is a list of
$X$-cobordisms $W_{i}$ with a given parenthesation (that is, a
word in the $X$-cobordisms $W_{i}$) such that
\begin{equation*}
  W=W_{1}\circ\cdots\circ W_{n}
\end{equation*}
with parenthesis in the middle as prescribed by $l$. The morphisms
$(W,l)\longrightarrow (W',l')$ are all $X$-diffeomorphisms
$W\longrightarrow W'$. Note that that we don't require that the
morphisms preserve the decomposition of $W$ into $W_{i}$,
nevertheless the associator $X$-diffeomorphisms $\Delta_{W,V,U}$
\emph{do preserve it} so we are in conditions to reuse the
argument of the coherence theorem and conclude that $\mathcal{G}$
is thin.

Now apply the above quotient construction to get a well defined
category (that will remain nameless). Note that any HQFT
trivially satisfies the conditions of the universal property
lemma, and thus, by a similar argument on the lemma it has a well
defined value on these quotient objects -- it does not define a
functor since an HQFT does not assign values to
$X$-diffeomorphisms between $X$-cobordisms. It is also easy to
see that the boundary of an $X$-diffeomorphism continues to be
well-defined.

From now on whenever we talk of $X$-cobordisms or
$X$-diffeomorphisms between them we shall always mean the objects
and morphisms in this quotient category.

Lest the reader be worried that we have introduced a whole series
of copies of a unique $X$-cobordism we hasten to add that we have
done so but in the end (in the last step, more precisely) these
will all be identified with one another.

\smallskip
\textbf{Step 4: Joining diffeomorphisms and cobordisms.}

Consider the wu-category $\Catname{C}(n,X)$ having for objects
$X$-manifolds of dimension $n$ and morphisms $X$-cobordisms. By
axiom \ref{axm:ax3} with $\psi$ the identity, and the above step,
given an HQFT $\tau$ there is a unique functor
$\Catname{C}(n,X)\longrightarrow\Catname{Mod}(R)$. It is also easy
to see that a map $\tau\longrightarrow\theta$ between HQFT's
induces a unique natural transformations between the corresponding
functors $\Catname{C}(n,X)\longrightarrow\Catname{Mod}(R)$. If we
consider the discrete category $\Catname{Man}(n,X)$ of $n$
dimensional $X$-manifolds there are functors
$\Catname{Man}(n,X)\longrightarrow\Catname{Diff}(n,X)$ and
$\Catname{Man}(n,X)\longrightarrow\Catname{C}(n,X)$, this one
taking the identities to the identity cobordisms, so that we can
form the push-out square in figure \ref{diag:pushoutsquare} to
obtain a wu-category.
\begin{figure}[hbtp]
  $\vcenter{\xymatrix{
    \Catname{Man}(n,X) \ar[d] \ar[r] & \Catname{Diff}(n,X) \ar[d] \\
    \Catname{C}(n,X) \ar[r] & \Catname{C}(n,X)\coprod_{\Catname{Man}(n,X)}\Catname{Diff}(n,X)   }}$
  \caption{}
  \label{diag:pushoutsquare}
\end{figure}

\smallskip
Push-outs are quite delicate objects and may or may not exist in
any given category, so we have to justify the use of the term. But
in this case, since $\Catname{Man}(n,X)$ sits very neatly inside
the other two categories we can give a very concrete recipe to
construct this particular push-out.

Consider the graph where the vertices are $X$-cobordisms and the
edges are $X$-diffeomorphisms \emph{plus} the $X$-cobordisms and
take the free wu-category on this graph. The objects of this
category are the vertices of the graph, i.e. $X$-manifolds, but
the morphisms are strings of $X$-cobordisms and
$X$-diffeomorphisms, composition being given by concatenation of
strings.

Denote by $<\!\psi\!>$ an $X$-diffeomorphism and by $<\!(W,F)\!>$
an $X$-cobordism in this free wu-category. The required push-out
is now obtained by taking the quotient by the congruence generated
by the relation
\begin{equation*}
\begin{split}
  <\!\psi\!><\!\phi\!> &\sim<\!\psi\phi\!> \\
  <\!(W,F)\!><\!(V,G)\!> &\sim<\!(W,F)\circ(V,G)\!> \\
  <\!(I\times M,1_{g})\!> &\sim<\! 1_{(M,g)}\!>
\end{split}
\end{equation*}

Because of these relations there are inclusion functors
\begin{align*}
  i:\Catname{C}(n,X) &\longrightarrow\Catname{C}(n,X)\coprod_{\Catname{Man}(n,X)}\Catname{Diff}(n,X) \\
  j:\Catname{Diff}(n,X) &\longrightarrow\Catname{C}(n,X)\coprod_{\Catname{Man}(n,X)}\Catname{Diff}(n,X)
\end{align*}

The reader is invited to check that these satisfy the universal
property of push-outs detailed in the diagram of figure
\ref{diag:pushoutunivprop}, where the push-out object is denoted
simply by $\mathcal{A}$.
\begin{figure}[hbtp]
  $\vcenter{\xymatrix{
  \Catname{Man}(n,X) \ar[d] \ar[r] & \Catname{Diff}(n,X) \ar[d]_{j} \ar@/^/[ddr] &  \\
  \Catname{C}(n,X) \ar[r]^{i} \ar@/_/[drr] & \mathcal{A} \ar@{-->}[dr] & \\
                &               & \mathcal{B}    }}$
  \caption{}
  \label{diag:pushoutunivprop}
\end{figure}

Since by axioms \ref{axm:ax1} and \ref{axm:ax2} there is a
symmetric monoidal functor
\begin{equation*}
  \Catname{Diff}(n,X)\longrightarrow\Catname{Mod}(R)
\end{equation*}
by the universal property of push-outs there is a unique functor
\begin{equation*}
  \Catname{C}(n,X)\coprod_{\Catname{Man}(n,X)}\Catname{Diff}(n,X)\longrightarrow\Catname{Mod}(R)
\end{equation*}

In what follows the remaining properties desirable for a category
of homotopy cobordisms will be given by quotienting this category
by congruences, or relations on strings of $X$-diffeomorphisms and
$X$-cobordisms. In other words, the subsequent steps will output a
category where the morphisms \emph{are} alternating strings of
$X$-diffeomorphisms and $X$-cobordisms, just with some added
relations and this fully justifies the fact that in the proofs of
the first section we have considered these two cases one at a
time.

\smallskip
\textbf{Step 5: Monoidal structure.}

Now for the monoidal structure. Axiom \ref{axm:ax5} tell us that
$X$-cobordisms should be considered on an equal footing with
$X$-diffeomorphisms as regards the monoidal structure. Since the
$X$-diffeomorphisms already have a monoidal structure what we will
do is make the associator $\alpha_{M,N,P}$, the unital structural
maps $l_M$ and $r_M$ and the symmetry $\sigma_{M,N}$, the
associator, the unital maps and the symmetry also for the
$X$-cobordisms. This amounts to imposing the naturality
constraints. So take the quotient by the congruence generated
by\footnote{From now on we omit the braces $<\!>$ in the
morphisms. Composition of an $X$-diffeomorphism and an
$X$-cobordism is denoted by concatenation.}
\begin{equation*}
\begin{split}
 \alpha_{M,N,P}(W\amalg (V\amalg U)) &\sim((W\amalg V)\amalg U)\alpha_{M',N',P'} \\
 l_M W &\sim (\emptyset\amalg W)l_N \\
 r_M W &\sim (W\amalg\emptyset)r_N \\
 \sigma_{M,N}(V\amalg W) &\sim (W\amalg V)\sigma_{M',N'}
\end{split}
\end{equation*}

Axiom \ref{axm:ax5} and Diagram \ref{diag:cobordassoc} in step 2
(and the comments following it) show that for a given HQFT $\tau$
there is a unique symmetric monoidal functor from this (nameless)
quotient into $\Catname{Mod}(R)$.

\smallskip
\textbf{Step 6: Relations.}

First let me define the notion of a monoidal congruence in a
category, the means to define quotients of monoidal categories.
For the definition of a congruence in a category the reader can do
no better than to turn to the book \cite{MAC}. A monoidal
congruence is a congruence with the added property
\begin{equation*}
  f\sim f' \,\text{and}\, g\sim g'\Longrightarrow f\otimes g\sim
  f'\otimes g'
\end{equation*}

If $\sim$ is a monoidal congruence then we can define the quotient
$\mathcal{A}/\!\!\sim$, which is readily seen to be a monoidal
category, and there is a (strict) monoidal functor
$\mathcal{A}\longrightarrow\mathcal{A}/\!\!\sim$.

Monoidal congruences share all the properties of congruences. In
particular, the intersection of a family of monoidal congruences
is still a monoidal congruence and the identity congruence is
clearly monoidal. So it makes sense to speak of the monoidal
congruence generated by a relation.

Given this we finally (!) enforce axioms \ref{axm:ax3},
\ref{axm:ax4} and \ref{axm:ax6}, thereby giving birth to the
category $\Catname{HCobord}(n,X)$, by taking yet another quotient,
now by the monoidal congruence generated by the following
relation: For every $X$-diffeomorphism $\Psi$ between
$X$-cobordisms we have
\begin{equation*}
\Psi|_{\partial_{-}M}(V,G)\sim(W,F)\Psi|_{\partial_{+}W}
\end{equation*}

For $X$-cobordisms $(W,F)$ and $(V,G)$ glued along an
$X$-diffeomorphism $\psi$ we have
\begin{equation*}
  (W,F)\amalg_{\psi}(V,G)\sim(W,F)\psi(V,G)
\end{equation*}

And we also impose for the identity cobordisms
\begin{equation*}
\begin{split}
  V\circ(I\times M,1_g) &\sim V \\
  (I\times M,1_g)\circ W &\sim W \\
\end{split}
\end{equation*}

This group of relations is precisely the one listed in the
introduction, the relations needed to prove the propositions of
section 1.

The last relation makes of $(I\times M,1_g)$ an actual unit for
composition, and the first is the one that allows to kill the
added degeneracy introduced in step 3. Using the notation
introduced in that step, just take the ``identity''
$X$-diffeomorphism $1_{W}:(W,W_{i})\longrightarrow (W,W'_{j})$ and
apply the first relation. So the morphisms \emph{are} alternating
strings of (genuine) $X$-cobordisms and $X$-diffeomorphisms as
stated. It is also the first relation that identifies two distinct
$X$-cobordisms -- when they are diffeomorphic by a diffeomorphism
that is the identity on the boundary. In fact, by inspection of
the relations two $X$-cobordisms are equal iff they are
diffeomorphic by a diffeomorphism that is the identity on the
boundary. Therefore the proofs given in section
\ref{section:defsandgenprops} are completely justified.

\smallskip
To sum it all up we have

\begin{nntheorem}
There is a natural isomorphism of categories
\begin{equation*}
  \Catname{HQFT}(n,X)\cong[\Catname{HCobord}(n,X),\Catname{Mod}(R)]
\end{equation*}
\end{nntheorem}
\begin{proof}
The above step-by-step construction provides a map
\begin{equation*}
  \Catname{HQFT}(n,X)\longrightarrow[\Catname{HCobord}(n,X),\Catname{Mod}(R)]
\end{equation*}

It is not difficult to prove that it is indeed a functor and that
it has an inverse follows from the (trivial) reconstruction of an
HQFT from the respective functor
$\Catname{HCobord}(n,X)\longrightarrow\Catname{Mod}(R)$. The fact
that it is a natural isomorphism is also easy to prove and we
leave to the reader the task of providing the specific details.
\end{proof}

This theorem reduces the study of HQFT's to the study of symmetric
monoidal functors on $\Catname{HCobord}(n,X)$. In other words, if
the reader is wondering if the category $\Catname{HCobord}(n,X)$
is trivial or not, the above theorem says that it is as trivial as
the category of HQFT's.

\smallskip
\textbf{Acknowledgements:} Although the first page of this paper
bears my name this was by no means a solo work. It could not have
been completed without the help of, and the fruitful discussions
with, Professors Roger Picken and Louis Crane.

My warm thanks also to Marco Mackaay and Paul Turner whose
comments and suggestions helped me in polishing this paper.

\end{document}